\documentclass[11pt]{amsart}
\usepackage{tabularx,booktabs}
\usepackage{caption}
\usepackage{amsmath}
\usepackage{amsfonts}

\usepackage{amscd}
\usepackage{amsthm}
\usepackage{amssymb} \usepackage{latexsym}
\usepackage{eufrak}
\usepackage{euscript}
\usepackage{epsfig}
\usepackage{graphics}
\usepackage{array}
\usepackage{enumerate}
\usepackage{dsfont}
\usepackage{color}
\usepackage{wasysym}
\usepackage{hyperref}
\usepackage{pdfsync}

\newcommand{\bel}[1]{\begin{equation}\label{#1}}

\newcommand{\be}{\begin{equation}}

\newcommand{\ba}{\begin{eqnarray}}
\newcommand{\ea}{\end{eqnarray}}

\newcommand{\qe}{\end{equation}}
\newcommand{\R}{{\mathbb R}}

\newcommand{\Z}{{\mathbb Z}}

\newcommand{\NNG}{\mathcal{PC}_{\geq 0}}

\newcommand{\Hmm}[1]{\leavevmode{\marginpar{\tiny%
$\hbox to 0mm{\hspace*{-0.5mm}$\leftarrow$\hss}%
\vcenter{\vrule depth 0.1mm height 0.1mm width \the\marginparwidth}%
\hbox to
0mm{\hss$\rightarrow$\hspace*{-0.5mm}}$\\\relax\raggedright #1}}}

\newtheorem{theorem}{Theorem}[section]

\newtheorem{lemma}[theorem]{Lemma}
\newtheorem{corollary}[theorem]{Corollary}

\newtheorem{conjecture}[theorem]{Conjecture}
\newtheorem{prop}[theorem]{Proposition}

\begin{document}

\title[The first gap for total curvatures of planar graphs]{The first gap for total curvatures of planar graphs with nonnegative curvature}
\author{Bobo Hua}
\email{bobohua@fudan.edu.cn}
\address{School of Mathematical Sciences, LMNS, Fudan University, Shanghai 200433, China}

\author{Yanhui Su}
\email{suyh@fzu.edu.cn}
\address{College of Mathematics and Computer Science, Fuzhou University, Fuzhou 350116, China}

\begin{abstract} 
We prove that the total curvature of a planar graph with nonnegative combinatorial curvature is at least $\frac{1}{12}$ if it is positive.  Moreover, we classify the metric structures of ambient polygonal surfaces for planar graphs attaining this bound.
\end{abstract}
\maketitle
\tableofcontents
Mathematics Subject Classification 2010: 31C05, 05C10.


\par
\maketitle

\bigskip


\section{Introduction}\label{s:intro}
The combinatorial curvature for planar graphs, which stands as the generalized Gaussian curvature of piecewise flat manifolds, was introduced by many authors \cite{Neva70,Stone76,Gromov87,Ishida90}. Many interesting geometric and analytic results have been obtained since then, see e.g. \cite{Zuk97,Woess98,Higuchi01,BP01,HJL02,LPZ02,HigS03,SunYu04,RBK05,BP06,DeVosMohar07,ChenChen08,Zhang08,Chenbf09,Keller10,KP11,Keller11,Oh17}.

Let $G=(V,E,F)$ be a semiplanar graph embedded into a surface $S$ (without boundary) with the set of vertices $V,$ the set of edges $E$ and the set of faces $F,$ see \cite{HJL15} (It is called planar if $S$ is either the 2-sphere or the plane). We say that a semiplanar graph $G$ is a \emph{tessellation} of $S$ if the following hold, see e.g. \cite{Keller11}:
\begin{enumerate}[(i)]
\item Every face is homeomorphic to a disk whose boundary consists of finitely many edges of the graph.
\item Every edge is contained in exactly two different faces.
\item For any two faces whose closures have non-empty intersection, the intersection is either a vertex or an edge.
\end{enumerate} We always consider semiplanar graphs which are tessellations of surfaces. 
For a semiplanar graph $G$, the \emph{combinatorial curvature} at the
vertex is defined as
\begin{equation}\label{def:comb}\Phi(x)=1-\frac{\deg(x)}{2}+\sum_{\sigma\in F:x\in \overline{\sigma}}\frac{1}{\deg(\sigma)},\quad x\in V,\end{equation} where the summation is taken over all faces $\sigma$ whose closure $\overline{\sigma}$ contains $x$ and $\deg(\cdot)$ denotes the degree of a vertex or a face, see Section~\ref{s:pre}. There is a natural metric space associated to the semiplanar graph $G$ embedded into $S,$ called the \emph{polygonal surface} of $G$ and denoted by $S(G)$: Replace each face of $G$ by a regular Euclidean polygon of side length one with same facial degree, glue these polygons along the common edges and consider the induced metric on $S.$ Note that $S(G)$ is piecewise flat and hence one can define the generalized Gaussian curvature $K$ of $S(G)$ which is concentrated on the vertices. It is well-known that for any $x\in V,$ \begin{equation}\label{eq:comgau}\Phi(x)=\frac{1}{2\pi}K(x),\end{equation} where $K(x)$ is the generalized Gaussian curvature (or the angle defect) at $x,$ i.e. the difference of 
$2\pi$ and the total angle at $x$ in the metric space $S(G),$ see e.g. \cite{Alex05}. In this paper, we only consider the combinatorial curvature of semiplanar graphs and simply call it the curvature if it is clear in the context.

For a smooth surface with absolutely integrable Gaussian curvature, its total curvature encodes the global geometric information of the space, see \cite{SST03}. 
For a semiplanar graph $G,$ we denote by $$\Phi(G)=\sum_{x\in V}\Phi(x)$$ the total curvature of $G$ whenever the summation converges absolutely.  In case of finite graphs, the Gauss-Bonnet theorem reads as, see e.g. \cite{DeVosMohar07},
\begin{equation}\label{eq:GBfinite}\Phi(G)=\chi(S(G)),\end{equation} where $\chi(\cdot)$ denotes the Euler characteristic of a surface.
For an infinite semiplanar graph $G,$ the Cohn-Vossen type theorem, see \cite[Theorem~1.3]{DeVosMohar07} or \cite[Theorem~1.6]{ChenChen08}, yields that
\begin{equation}\label{CVthm}\Phi(G)\leq\chi(S(G)),\end{equation}
whenever $\sum_{x\in V}\min\{\Phi(x),0\}$ converges.

In this paper, we study total curvatures of semiplanar graphs with nonnegative combinatorial curvature. Note that a semiplanar graph $G$ has nonnegative combinatorial curvature if and only if the polygonal surface $S(G)$ is a generalized convex surface, see \cite{BuragoGromovPerelman92,BuragoBuragoIvanov01,HJL15}. We denote by $\NNG$ the set of infinite semiplanar graphs with nonnegative combinatorial curvature. By the Cohn-Vossen type theorem, the inequality \eqref{CVthm} yields that $0\leq \Phi(G)\leq 1,$ for any $G\in \NNG.$ 
Due to the discrete nature, R\'eti proposed the following conjecture, see \cite[Conjecture~2.1]{HuaLin16}: $$\tau_1:=\inf\left\{\Phi(G):G\in \NNG,\Phi(G)>0\right\}=\frac16$$ and the minimum is attained by the graph consisting of a pentagon and infinitely many hexagons. We call $\tau_1$ the first gap of the total curvature for planar graphs with nonnegative curvature.


For a semiplanar graph $G=(V,E,F)$ with nonnegative combinatorial curvature, we denote by $$T(G):=\{x\in V:\Phi(x)\neq 0\}$$ the set of vertices with non-vanishing curvature. Some crucial geometric information is contained in the structure of $T(G).$ Chen and Chen \cite{ChenChen08,Chenbf09} obtained
an interesting result that the curvature vanishes outside a finite subset of vertices in an infinite semiplanar graph with nonnegative combinatorial curvature. 
\begin{theorem}[Theorem~1.4 in \cite{ChenChen08}, Theorem~3.5 in \cite{Chenbf09}]\label{thm:ccthm} For a semiplanar graph $G$ with nonnegative combinatorial curvature, $T(G)$ is a finite set.
\end{theorem}

In the previous paper \cite{HuaSu17a}, we determine all possible values of total curvature of semiplanar graphs with nonnegative curvature. 
\begin{theorem}[Theorem~1.1 in \cite{HuaSu17a}]\label{mainthm2}
The set of all values of total curvatures of infinite semiplanar graphs with nonnegative combinatorial curvature is given by $$\left\{\frac{i}{12}:\ 0\leq i\leq 12, i\in \Z\right\}.$$
\end{theorem} As a corollary, one immediately obtains $\tau_1=\frac{1}{12}$ which answers R\'eti's question. The proof of the theorem crucially uses Theorem~\ref{thm:ccthm} and the Gauss-Bonnet theorem on compact subsets with boundary.
However, it leaves the structure of the graph $G$ or the subset $T(G)$ in a black box, that is, we don't know any further information of the structures of semiplanar graphs attaining the first gap $\tau_1$ by Theorem~\ref{mainthm2}.

\begin{figure}[htbp]
\begin{center}
\includegraphics[width=0.9\linewidth]{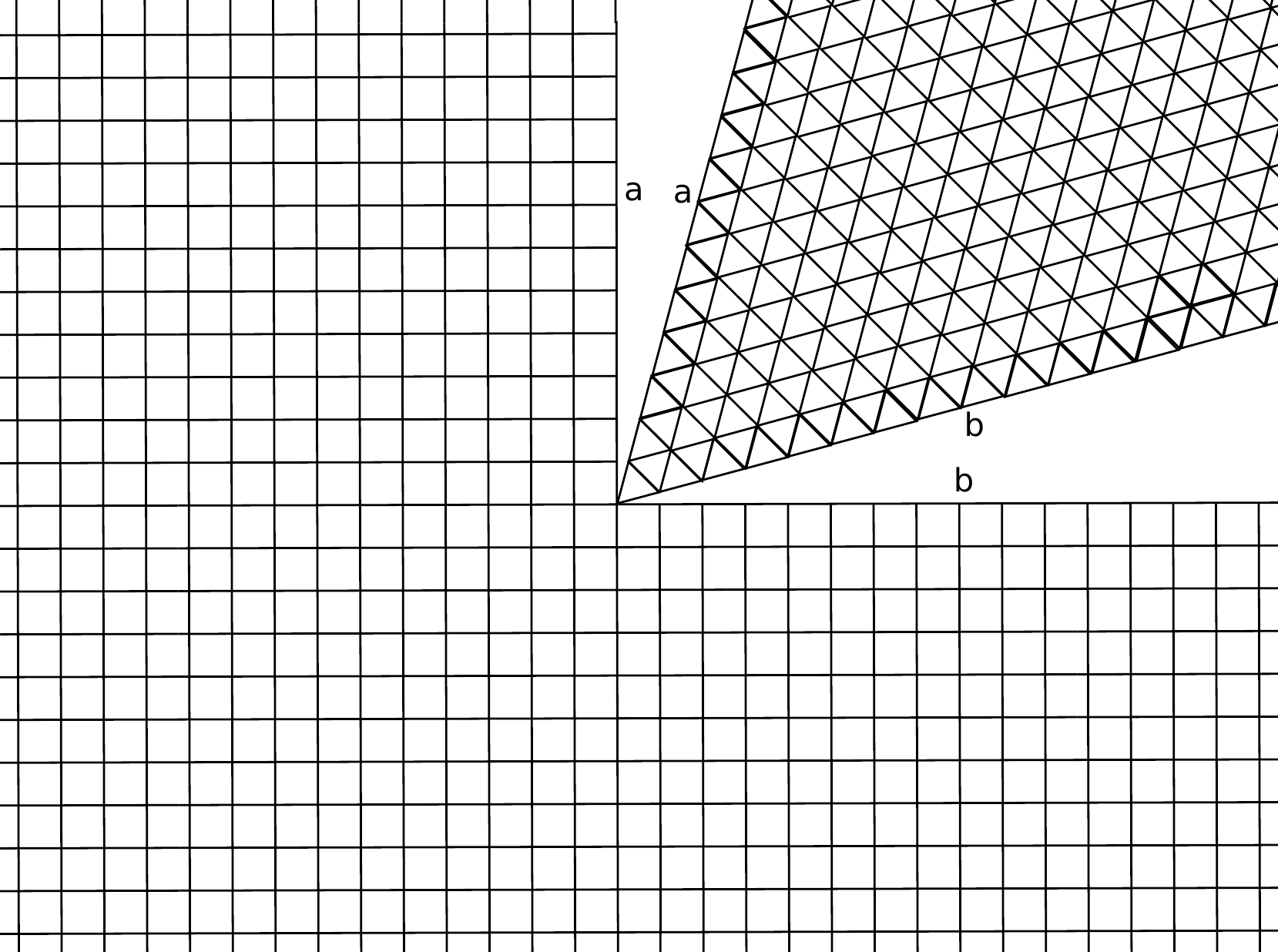}
\caption{\small A graph with a single vertex of curvature $\frac{1}{12}.$}
\label{fig2}
\end{center}
\end{figure}

\begin{figure}[htbp]
\begin{center}
\includegraphics[width=0.9\linewidth]{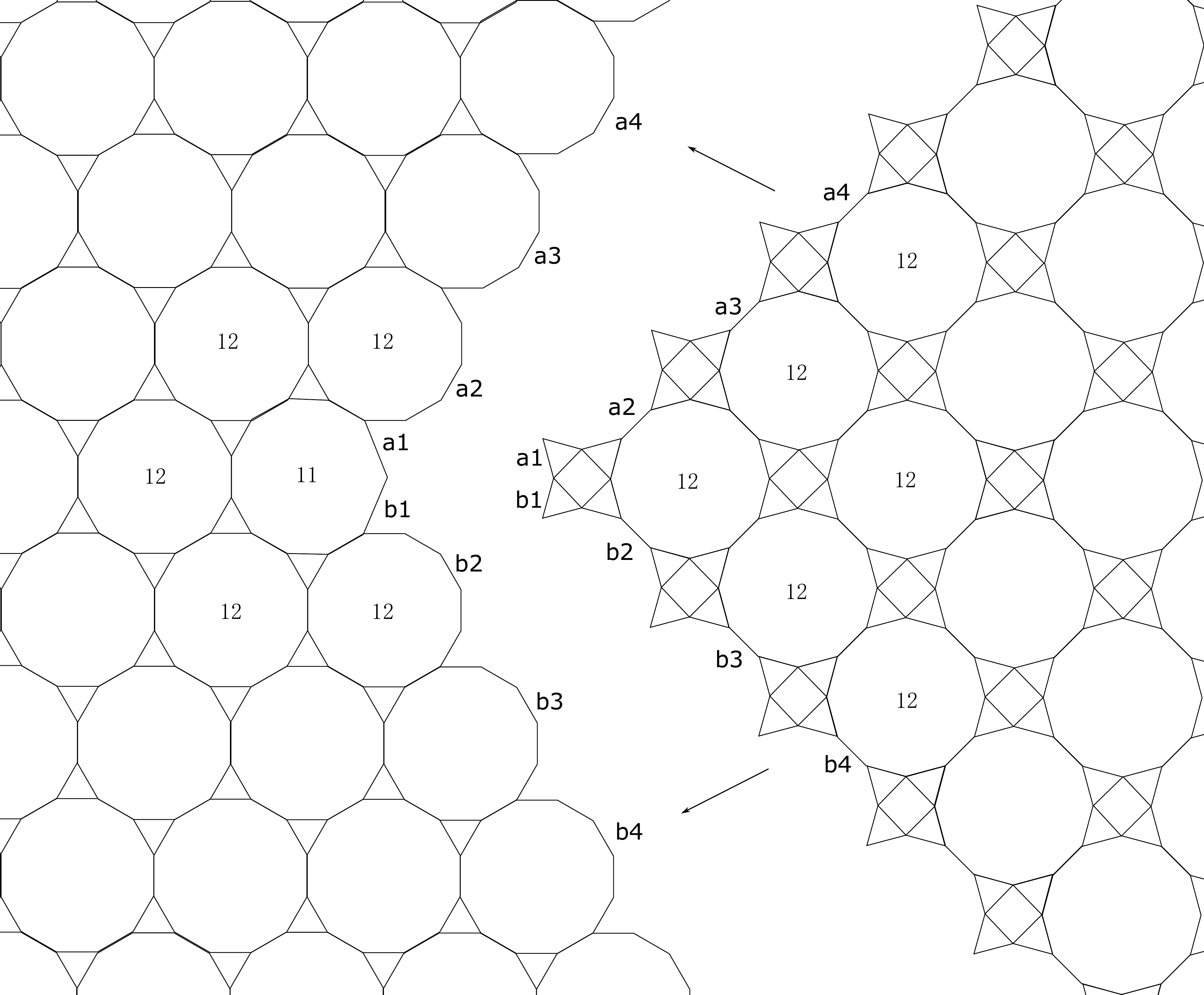}
\caption{\small A graph with eleven vertices (the vertices of the hendecagon) of curvature $\frac{1}{132}.$}
\label{fig5}
\end{center}
\end{figure}

In this paper, we adopt a different approach to study the first gap of the total curvature and give the answer to R\'eti's question independent of Theorem~\ref{mainthm2}, i.e. it does not use Theorem~\ref{thm:ccthm} and the Gauss-Bonnet theorem. On one hand, we show that the total curvature in the class of $\NNG$ is at least $\frac{1}{12}$ if it is positive. On the other hand, we construct several examples whose total curvature are $\frac{1}{12},$ see Figure~\ref{fig2} and Figure~\ref{fig5} (see Section~\ref{s:examples} for more examples). For any figure in this paper, all vertices, edges or sides with same labeling are identified with each other. The proof strategy is straight-forward and involves tedious case studies, see Section~\ref{s:thm1}: For a vertex with non-vanishing curvature, if the curvature of the vertex is less than $\frac{1}{12},$ then we try to find some nearby vertices with non-vanishing curvature such that the sum of these curvatures is at least $\frac{1}{12}$ and prove the results case by case. This approach provides us more information of the structures of semiplanar graphs attaining the first gap of the total curvature. In particular, we get the classification of metric structures for polygonal surfaces of such semiplanar graphs.

\begin{figure}[tb]
\begin{minipage}[b]{\textwidth}
\centering
\includegraphics[width=11cm,height=3.8cm]{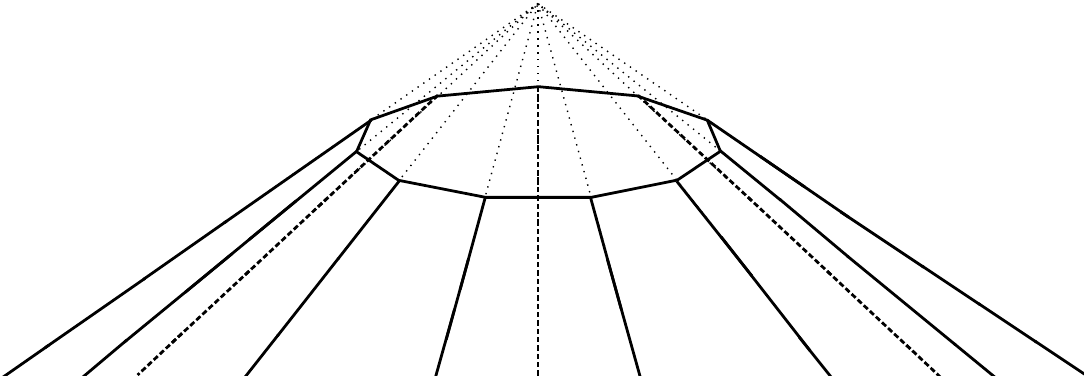}
\end{minipage}
\caption{\small A ``frustum" with a hendecagon base.
}
 \label{fig-11}
\end{figure}

\begin{theorem}\label{mainthm} 
$$\tau_1=\frac{1}{12}.$$
A semiplanar graph $G$ with nonnegative curvature satisfies $\Phi(G)=\frac{1}{12}$ if and only if the polygonal surface $S(G)$ is isometric to either
\begin{enumerate}[(a)]
\item a cone with the apex angle $\theta=2\arcsin\frac{11}{12},$ or
\item a ``frustum" with a hendecagon base, see Figure \ref{fig-11}.
\end{enumerate}
\end{theorem}

Note that in case of $(a)$ in the above theorem the graph has a single vertex at the apex of the polygonal surface with the curvature $\frac{1}{12},$ see e.g. Figure~\ref{fig2}; In case of $(b)$ it has eleven vertices, on the boundary of the base of the polygonal surface, with the curvature $\frac{1}{132}$ for each, see e.g. Figure~\ref{fig5}. It is usually hard to classify the graph/tessellation structures even in the class of planar tilings with vanishing curvature, see  \cite{GrunbaumS87}.
However, our result indicates that semiplanar graphs attaining the first gap of the total curvature have rigid metric structure for ambient polygonal surfaces. 

In the last section, we include some applications of the main result.
On one hand, we prove that for an infinite semiplanar graph $G$ with nonnegative curvature, the induced subgraph on $T(G)$ has at most $14$ connected components, see Corollary~\ref{coro:concomp}. An example with $12$ connected components is constructed in Figure~\ref{12-component}.
On the other hand, we obtain the first gap of the total curvature for semiplanar graphs with boundary and with nonnegative curvature is $\frac{1}{12},$ see \cite[Theorem~5.5]{HuaSu17a} for an alternative proof. For precise definitions and the terminology used, we refer to Subsection~\ref{s:boundary}. Let $G$ be a semiplanar graph with boundary and with nonnegative curvature, and $S(G)$ the polygonal surface of $G.$ Consider the doubling constructions of $S(G)$ and $G,$ see e.g. \cite[Section~5]{Perelman91}. Let $\widetilde{S(G)}$ be the
double of ${S}(G)$, that is, $\widetilde{S(G)}$ consists of two copies of ${S}(G)$ glued along the boundaries via the identity map restricted on $\partial S(G).$
This induces the doubling graph of $G,$ denoted by $\widetilde{G}.$ By the definition of the curvature for semiplanar graphs with boundary, one can show that $\widetilde{G}$ has nonnegative curvature and $\Phi(\widetilde{G})=2\Phi(G).$ This yields that the total curvature of $G$ is an integral multiple of $\frac{1}{24}$ by Theorem~\ref{mainthm2}, which suggests that the first gap for semiplanar graphs with boundary could be $\frac{1}{24}.$ However, using the rigidity of metric structures for seimplanar graphs without boundary, Theorem~\ref{mainthm}, we can conclude that $\Phi(\widetilde{G})\neq\frac{1}{12}$ and hence $\Phi(G)\geq \frac1{12}$ by Theorem~\ref{mainthm2}, see Corollary~\ref{coro:withboundary}.

The paper is organized as follows: In next section, we recall some basics on the combinatorial curvature for semiplanar graphs.
Section~\ref{s:thm1} is devoted to the proof of Theorem~\ref{mainthm}.
In Section~\ref{s:examples}, we give some examples to show the sharpness of Theorem~\ref{mainthm}. In the last section, we give some applications of Theorem~\ref{mainthm}.


\section{Preliminaries}\label{s:pre}
Let $G=(V,E,F)$ be a semiplanar graph induced by an embedding of a graph $(V,E)$ into a (possibly un-orientable) surface $S$ without boundary. Any connected component of the complement of the embedding image of the graph $(V,E)$ into $S$ is called a face. We only consider the appropriate embedding such that $G$ is a tessellation of $S,$ see the definition in the introduction.

We say that a vertex is incident to an edge (similarly, an edge is incident to a face, or a vertex is incident to a face) if the former is a subset of the closure of the latter. Two vertices are called neighbors if there is an edge connecting them. We denote by $\deg(x)$ the number of neighbors of a vertex $x,$ and by $\deg(\sigma)$ the number of edges incident to a face $\sigma$ (equivalently, the number of vertices incident to $\sigma$). For a tessellation, we always assume that $3\leq \deg(x)<\infty$ and $3\leq\mathrm{deg}(\sigma)<\infty$ for any vertex $x$ and face $\sigma.$
 Two edges (two faces resp.) are called adjacent if there is a vertex (an edge resp.) incident to both of them.
The combinatorial distance between two vertices $x$ and $y,$ denote by $d(x,y),$ is defined as the minimal length of walks from $x$ to $y,$ i.e. the minimal number $n$ such that there is $\{x_{i}\}_{i=1}^{n-1}\subset V$ satisfying $x\sim x_1\sim\cdots\sim x_{n-1}\sim y.$ We denote by $B_r(x):=\{y\in V: d(y,x)\leq r\},$ $r\geq 0,$ the ball of radius $r$ centered at the vertex $x.$

Given a semiplanar graph $G=(V,E,F)$ embedded into a surface $S,$ it associates with a unique metric space $S(G),$ called the polygonal surface defined in the introduction. The combinatorial curvature at a vertex is defined to be proportional to the generalized Gaussian curvature on $S(G),$ see \eqref{eq:comgau}.
In this paper, we study total curvatures of infinite semiplanar graphs with nonnegative combinatorial curvature. For our purposes, it suffices to consider those with positive total curvature (otherwise the total curvature vanishes). By \cite[Theorem 3.10]{HJL15}, these graphs are planar, namely, the ambient spaces are homeomorphic to $\R^2.$

In a semiplanar graph, a pattern of a vertex $x$ is defined as a vector $(\mathrm{deg}(\sigma_1),\mathrm{deg}(\sigma_2),\cdots,\mathrm{deg}(\sigma_{N})),$ where $N=\deg(x),$ $\{\sigma_i\}_{i=1}^N$ are the faces which $x$ is incident to, and $\mathrm{deg}(\sigma_1)\leq\mathrm{deg}(\sigma_2)\leq\cdots\leq\mathrm{deg}(\sigma_{N}).$ For simplicity, we always write $$x=(\mathrm{deg}(\sigma_1),\mathrm{deg}(\sigma_2),\cdots,\mathrm{deg}(\sigma_{N}))$$ to indicate the pattern of the vertex $x.$

Table \ref{tabl1} is the list of all
possible patterns of a vertex with positive curvature (see
\cite{DeVosMohar07,ChenChen08}); Table \ref{tabl2} is the list of all
possible patterns of a vertex with vanishing curvature (see \cite{GrunbaumS87,ChenChen08}).

\begin{table}
\refstepcounter{table}\label{tabl1}
\begin{tabular}{|lc|lr|}
\hline
Patterns &&$\Phi(x)$&\\
    \hline
 $(3,3,k)$ & $3\leq k$&$1/6+1/k$&\\
 $(3,4,k)$  & $4\leq k$  &$1/12+1/k$&\\
 $(3,5,k)$ & $5\leq k$&$1/30+1/k$&\\
$(3,6,k)$&$6\leq k$&$1/k$&\\
$(3,7,k)$&$7\leq k\leq41$&$1/k-1/42$&\\
$(3,8,k)$&$8\leq k\leq 23$&$1/k-1/24$&\\
$(3,9,k)$&$9\leq k\leq 17$&$1/k-1/18$&\\
$(3,10,k)$&$10\leq k\leq 14$&$1/k-1/15$&\\
$(3,11,k)$&$11\leq k\leq 13$&$1/k-5/66$&\\
$(4,4,k)$&$4\leq k$&$1/k$&\\
$(4,5,k)$&$5\leq k\leq 19$&$1/k-1/20$&\\
$(4,6,k)$&$6\leq k\leq 11$&$1/k-1/12$&\\
$(4,7,k)$&$7\leq k\leq 9$&$1/k-3/28$&\\
$(5,5,k)$&$5\leq k\leq 9$&$1/k-1/10$&\\
$(5,6,k)$&$6\leq k\leq 7$&$1/k-2/15$&\\
$(3,3,3,k)$&$3\leq k$&$1/k$&\\
$(3,3,4,k)$&$4\leq k\leq 11$&$1/k-1/12$&\\
$(3,3,5,k)$&$5\leq k\leq 7$&$1/k-2/15$&\\
$(3,4,4,k)$&$4\leq k\leq 5$&$1/k-1/6$&\\
$(3,3,3,3,k)$&$3\leq k\leq5$&$1/k-1/6$&\\
\hline
\multicolumn{4}{c}{}\\    
\end{tabular}

\textbf{\tablename~\thetable.} The patterns of a vertex with positive curvature.
\end{table}

\begin{table}
\refstepcounter{table}\label{tabl2}
\begin{tabular}{lllll}
\hline
$(3,7,42),$ & $(3,8,24),$ & $(3,9,18),$ & $(3,10,15),$ & $(3,12,12),$\\
$(4,5,20),$ & $(4,6,12),$ & $(4,8,8),$ & $(5,5,10),$ & $(6,6,6),$\\
$(3,3,4,12),$ & $(3,3,6,6),$& $(3,4,4,6),$ & $(4,4,4,4),$ &$(3,3,3,3,6),$\\
$(3,3,3,4,4),$&$(3,3,3,3,3,3).$&&&\\
\hline
\multicolumn{1}{c}{}\\    
\end{tabular}

\textbf{\tablename~\thetable.} The patterns of a vertex with vanishing curvature.
\end{table}

%
%

The following lemma is useful in this paper, see \cite[Lemma~2.5]{ChenChen08}.
\begin{lemma}\label{lemma}
If there is a face $\sigma$ such that $\mathrm{deg}(\sigma)\geq43$ and $\Phi(x)\geq0$ for
any vertex $x$ incident to $\sigma,$ then
$$\sum_{x}\Phi(x)\geq1,$$ where the summation is taken over all vertices $x$ incident to $\sigma.$
\end{lemma}

In fact, for an infinite semiplanar graph $G$ with nonnegative curvature if there is a face whose degree is at least $43,$ then the graph has rather special
structure, see \cite[Theorem~2.10]{HJL15}. In particular, the total curvature of $G$ is equal to $1$ and the polygonal surface $S(G)$ is isometric to a half flat-cylinder
in $\mathbb{R}^3.$

\section{Proof of Theorem \ref{mainthm}}\label{s:thm1}
In this section, we give the proof of the Theorem \ref{mainthm}. Let $G=(V,E,F)$ be an infinite planar graph with nonnegative curvature and positive total curvature. The main strategy is to prove the results case by case: Pick a vertex $A$ with non-vanishing curvature and find many nearby vertices with non-vanishing curvature such that the sum of their curvatures is at least $\frac{1}{12}$. 

\begin{proof}[Proof of Theorem \ref{mainthm}]
By Lemma \ref{lemma} and \eqref{CVthm}, or \cite[Theorem~2.10]{HJL15}, it suffices to consider the case that 
$$\mathrm{deg}(\sigma)\leq 42,\quad \forall\ \sigma \in F.$$ 

For $A\in V$ with $\Phi(A)\neq 0,$ we consider all possible patterns of $A$ as follows (other vertices are labelled as in Figures):
\begin{description}
\item[Case 1] $A=(3,3,k)$. In this case,  $\Phi(A)=\frac{1}{6}+\frac{1}{k}>\frac{1}{12}$.
\item[Case 2]  $A=(3,4,k)$. In this case,  $\Phi(A)=\frac{1}{12}+\frac{1}{k}>\frac{1}{12}$.
\item [Case 3] $A=(3,5,k)$. In this case,  $\Phi(A)=\frac{1}{30}+\frac{1}{k}$. If $k\leq 19$, then $\Phi(A)>\frac{1}{12}$.
So it suffices to consider $k\geq 20.$ We denote by $B$ the neighbor of $A$ which is incident to the triangle and the $k$-gon, from this case to Case 9 in the below.
The possible pattern of $B$ is $(3,5,k)$, $(3, 6,k),$ $(3, 7,k), k\leq 42,$ $(3,8,k),k\leq 24,$ 
or $(3,3,3,k)$.

\begin{figure}[tb]
\begin{minipage}[b]{0.4\textwidth}
\centering
\includegraphics[width=4.2cm,height=4.2cm]{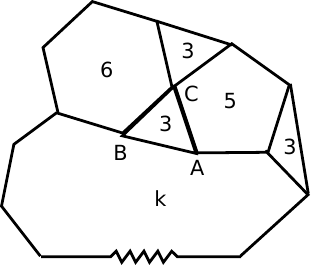}
(a)
\end{minipage}
\begin{minipage}[b]{0.4\textwidth}
\centering
\includegraphics[width=4.2cm,height=4.2cm]{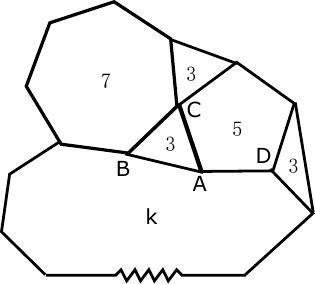}
(b)
\end{minipage}
\begin{minipage}[b]{0.4\textwidth}
\centering
\includegraphics[width=4.2cm,height=4.2cm]{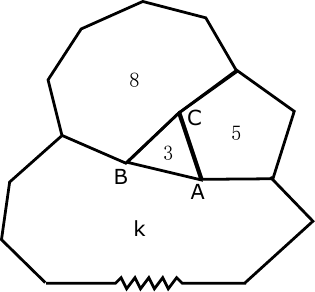}
(c)
\end{minipage}
\begin{minipage}[b]{0.4\textwidth}
\centering
\includegraphics[width=4.2cm,height=4.2cm]{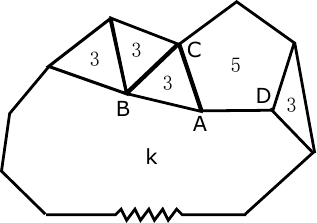}
(d)
\end{minipage}
\caption{\small Case 3.}
 \label{fig-case3}
\end{figure}

\begin{itemize}
\item If $B=(3,5,k)$, then $\Phi(B)=\frac{1}{30}+\frac{1}{k}$. For
$k\leq 42,$ $\Phi(A)+\Phi(B)>\frac{1}{12}$.

\item If $B=(3,6,k),$ then the only nontrivial case is shown in Figure \ref{fig-case3}(a). In this case, $\Phi(B)=\frac{1}{k}$. Since 
$C=(3,3,5,6) $, $\Phi(C)=\frac{1}{30}$ and $k\leq 42$,  we have $\Phi(A)+\Phi(B)+\Phi(C)>\frac{1}{12}$.

\item If $B=(3,7,k),$ $k\leq 42$, then the only nontrivial case is shown in Figure \ref{fig-case3}(b).
In this case $\Phi(B)=\frac{1}{k}-\frac{1}{42}$, $\Phi(C)=\frac{1}{105}$. If $k\leq20$, then $\Phi(A)+\Phi(B)>\frac{1}{12}$. If $k>20$, then $D$ has to be $(3,5,k)$, 
 $\Phi(D)=\frac{1}{30}+\frac{1}{k}$, which yields that $\Phi(A)+\Phi(D)>\frac{1}{12}$.

\item If $B=(3,8,k),$ $k\leq 24$, see Figure \ref{fig-case3}(c), then $\Phi(C)>\frac{1}{12}$ which implies that $\Phi(A)+\Phi(C)>\frac{1}{12}$.

\item If  $B=(3,3,3,k),$ then the only nontrivial case is shown in Figure \ref{fig-case3}(d). In this case $\Phi(B)=\frac{1}{k}$.  If $k\leq20$, then $\Phi(A)+\Phi(B)>\frac{1}{12}$. If $k>20$, then the pattern of $D$ has to be $(3,5,k)$ and 
 $\Phi(D)=\frac{1}{30}+\frac{1}{k}.$ Hence
$\Phi(A)+\Phi(D)>\frac{1}{12}$.
\end{itemize}

\item[Case 4] $A=(3,6,k)$. In this case,  $\Phi(A)=\frac{1}{k}$. For $k<12$, $\Phi(A)>\frac{1}{12}.$ For $k=12$, we can construct a graph with total curvature $\frac{1}{12},$ see Figure \ref{fig1}, whose polygonal surface is isometric to a cone with apex angle $2\arcsin\frac{11}{12}.$
So we only need to consider $k>12.$ The possible patterns of $B$ are
 $(3,6,k),$ $(3,7,k), k\leq 42,$ $(3,8,k), k\leq 24, $ $ (3,9,k), k\leq 18,$ $(3,10,k),k\leq 15,$
 $(3,11,k), k\leq13,$ and $(3,3,3,k)$.
 
 \begin{figure}[tb]
\begin{minipage}[b]{0.4\textwidth}
\centering
\includegraphics[width=4.2cm,height=4.2cm]{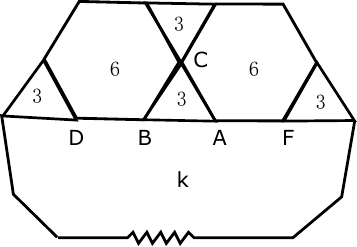}
(a)
\end{minipage}
\begin{minipage}[b]{0.4\textwidth}
\centering
\includegraphics[width=4.2cm,height=4.2cm]{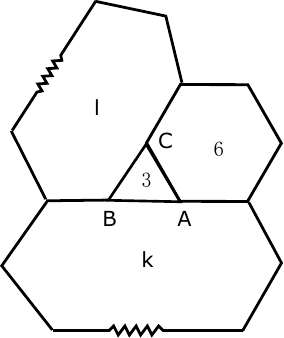}
(b)
\end{minipage}
\begin{minipage}[b]{0.4\textwidth}
\centering
\includegraphics[width=4.2cm,height=4.2cm]{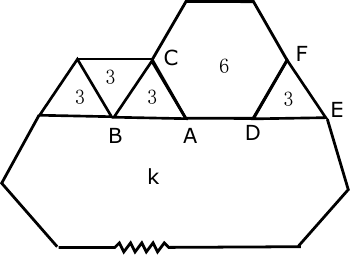}
(c)
\end{minipage}
\caption{\small Case 4.}
 \label{case4}
\end{figure}
 
 \begin{itemize}
\item
If $B=(3,6,k),$ then the only nontrivial case is shown in Figure \ref{case4}(a). In this case $\Phi(B)=\frac{1}{k}$ and $\Phi(C)=0.$
By the fact $k>12$, the patterns of $D$ and $F$ are $(3,6,k),$ which yields that $\Phi(D)=\Phi(F)=\frac{1}{k}$.
Hence $\Phi(A)+\Phi(B)+\Phi(D) +\Phi(F) =\frac{4}{k}>\frac{1}{12}$.
 
\item If $B=(3,l,k)$ for $l=7,8,9,10,11$, see Figure \ref{case4}(b), then $\Phi(C)=\frac{1}{l}>\frac{1}{12}$. This yields $\Phi(A)+\Phi(C)>\frac{1}{12}$.

\item If $B=(3,3,3,k),$ then the only nontrivial case is shown in Figure \ref{case4}(c). In this case $\Phi(B)=\frac{1}{k}$. Since $k>12$, $D=(3,6,k)$ and
 $\Phi(D)=\frac{1}{k}$.  Obviously, if $k\leq24$, then we have 
 $\Phi(A)+\Phi(B)+\Phi(D)>\frac{1}{12}$.
 If $k\geq25$, then the nontrivial patterns of $E$ are $(3,6,k)$, $(3,7,k)$ and $(3,3,3,k)$. If $E=(3,6,k)$ or $(3,3,3,k)$, then $\Phi(E)=\frac{1}{k}$,
thus $\Phi(A)+\Phi(B)+\Phi(D)+\Phi(E)=\frac{4}{k}>\frac{1}{12}$. If $E=(3,7,k)$, then $F=(3, 6,7)$, which implies $\Phi(F)=\frac{1}{7}$. Hence $\Phi(A)+\Phi(D)+\Phi(F)>\frac{1}{12}$.
\end{itemize}

\item[Case 5]  $A=(3,7,k), 7\leq k\leq41.$ In this case, $\Phi(A)=\frac{1}{k}-\frac{1}{42}$. If $k\leq 9$, then $\Phi(A)>\frac{1}{12}$.
So we only need to consider the case  $k\geq 10$. The possible patterns of  $B$ are
 $(3,7,k), k\leq 41,$ $(3,8,k), k\leq 24, $  $ (3,9,k), k\leq 18, $  $(3,10,k),k\leq 15,$
 $(3,11,k),k\leq 13, $ $(3,3,3,k),$ $(3,3,4,k),k\leq 12,$ and $(3,12,12).$

 \begin{figure}[tb]
\begin{minipage}[b]{0.4\textwidth}
\centering
\includegraphics[width=4.2cm,height=4.2cm]{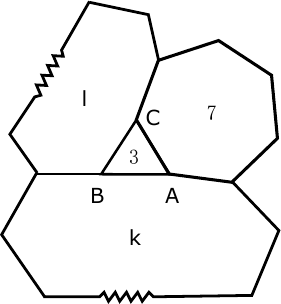}
(a)
\end{minipage}
\begin{minipage}[b]{0.4\textwidth}
\centering
\includegraphics[width=4.2cm,height=4.2cm]{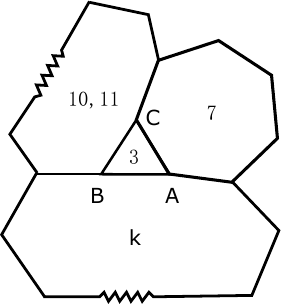}
(b)
\end{minipage}
\begin{minipage}[b]{0.4\textwidth}
\centering
\includegraphics[width=4.2cm,height=4.2cm]{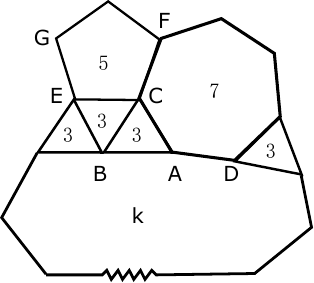}
(c)
\end{minipage}
\begin{minipage}[b]{0.4\textwidth}
\centering
\includegraphics[width=4.2cm,height=4.2cm]{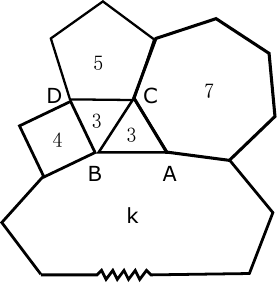}
(d)
\end{minipage}
\begin{minipage}[b]{0.4\textwidth}
\centering
\includegraphics[width=4.2cm,height=4.2cm]{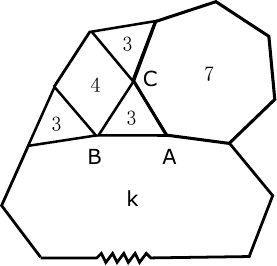}
(e)
\end{minipage}
\begin{minipage}[b]{0.4\textwidth}
\centering
\includegraphics[width=4.2cm,height=4.2cm]{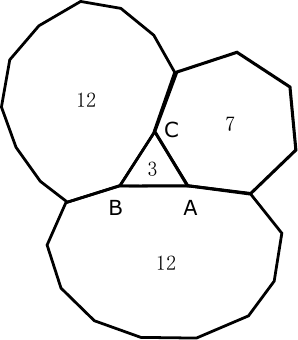}
(f)
\end{minipage}
\caption{\small Case 5.}
 \label{case5}
\end{figure}

\begin{itemize}
\item 
If $B=(3,l,k)$ for $l=7,8,9,$ see Figure \ref{case5}(a), then $\Phi(C)=\frac{1}{l}-\frac{1}{42}>\frac{1}{12}$. Hence $\Phi(A)+\Phi(C)>\frac{1}{12}$.


\item
If $B=(3,10,k)$ for $k\leq 15,$ see Figure \ref{case5}(b), then $\Phi(B)=\frac{1}{k}-\frac{1}{15}$ and $\Phi(C)=\frac{8}{105}$.
By  $k\leq 15$,  $\Phi(A)+\Phi(B)+\Phi(C)>\frac{1}{12}$.

\item
If $B=(3,11,k)$ for $k\leq 13,$ see also Figure \ref{case5}(b), then $\Phi(B)=\frac{1}{k}-\frac{5}{66}$ and $\Phi(C)=\frac{31}{462}$.
Since  $k\leq 13$, $\Phi(A)+\Phi(B)+\Phi(C)>\frac{1}{12}$.

\item
If $B=(3,3,3,k),$ see Figure \ref{case5}(c), then $\Phi(B)=\frac{1}{k}.$ 
The possible patterns of $C$ are $(3,3,3,7)$, $(3,3,4,7),$ and $(3,3,5,7)$. 
If $C=(3,3,3,7)$ or $(3,3,4,7)$, then $\Phi(C)\geq\frac{5}{84}$. So $\Phi(A)+\Phi(B)+\Phi(C)>\frac{1}{12}$.
If  $C=(3,3,5,7),$ then $\Phi(C)=\frac{1}{105 }$. For $k\leq 20$, $\Phi(A)+\Phi(B)+\Phi(C)>\frac{1}{12}$. So we only consider $k\geq21$. In this case, the pattern of $D$ has to be 
$(3,7,k)$ and $\Phi(D)=\frac{1}{k}-\frac{1}{42}.$ 
If $k\leq 24, $ then we have   $\Phi(A)+\Phi(B)+\Phi(C) +\Phi(D)>\frac{1}{12}$. If $k>24$, then the curvatures of other vertices on the $k$-polygon, except $B,A$ and $D,$ are at least $\frac{1}{k}-\frac{1}{42}$. Furthermore, we consider the vertex $E$. The possible patterns of $E$ are $(3,3,3,5)$, $(3,3,4,5)$, $(3,3,5,5)$, $(3,3,5,6)$ and $(3,3,5,7)$.
If $E=(3,3,3,5), $ $(3,3,4,5),$ $(3,3,5,5)$ or $(3,3,5,6),$ then
$\Phi(E)\geq\frac{1}{30}$ and $\sum\limits_{{x\in k\mbox{-}\mathrm{gon}}} \Phi(x)+\Phi(C)+\Phi(E)>\frac{1}{12},$ where the summation is taken over the vertices $x$ on the $k$-gon.
If  $E=(3,3,5,7), $  then
$\Phi(E)=\frac{1}{105}$ (in this case, $G$ is not a vertex of the $k$-gon, otherwise $G=(5,7,k)$ which has negative curvature). This yields that $\min\{\Phi(F), \Phi(G)\} \geq \frac{1}{105}$. Hence $\Phi(G)\geq\sum\limits_{{x\in k\mbox{-}\mathrm{gon}}} \Phi(x)+\Phi(C)+\Phi(E)+\Phi(F)+\Phi(G)>\frac{1}{12}.$

\item If $B=(3,3,4,k)$ for $k\leq 12$, then $\Phi(B)=\frac{1}{k}-\frac{1}{12}$. There are two subcases: The first one is shown in Figure \ref{case5}(d). 
Note that $k\leq 12$, if $C=(3,3,3,7)$ or $(3,3,4,7),$ then $\Phi(A)+\Phi(C)>\frac{1}{12}$.
So, the only nontrivial case is $C=(3,3,5,7)$ which yields $\Phi(D)\geq \frac{1}{30}$ and $\Phi(A)+\Phi(B)+\Phi(C)+\Phi(D)>\frac{1}{12}$;
 The other is depicted in Figure \ref{case5}(e). This yields that $\Phi(C)=\frac{5}{84}$ by $C=(3,3,4,7),$ and hence $\Phi(A)+\Phi(C)>\frac{1}{12}$.
 
\item If $B=(3,12,12)$, see Figure \ref{case5}(f), then $\Phi(A)=\frac{5}{84}$, $\Phi(B)=0$,  and $\Phi(C)=\frac{5}{84}$,
 which yields that $\Phi(A)+\Phi(C)>\frac{1}{12}$.
 \end{itemize}
 \item[Case 6] 
 $A=(3,8,k)$ for $8\leq k\leq23.$ In this case, $\Phi(A)=\frac{1}{k}-\frac{1}{24}$. The possible patterns of  $B$ are
  $(3,8,k), k\leq 23, $  $ (3,9,k), k\leq 18, $  $(3,10,k),k\leq 15,$
 $(3,11,k),k\leq 13, $ $(3,3,3,k),$ $(3,3,4,k),k\leq 12,$ and $(3,12,12).$
 
 \begin{figure}[tb]
\begin{minipage}[b]{0.4\textwidth}
\centering
\includegraphics[width=4.2cm,height=4.2cm]{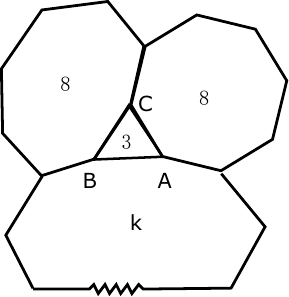}
(a)
\end{minipage}
\begin{minipage}[b]{0.4\textwidth}
\centering
\includegraphics[width=4.2cm,height=4.2cm]{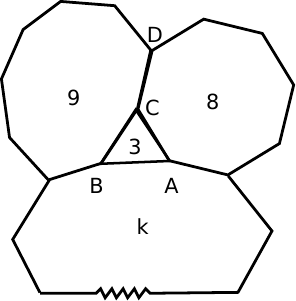}
(b)
\end{minipage}
\begin{minipage}[b]{0.4\textwidth}
\centering
\includegraphics[width=4.2cm,height=4.2cm]{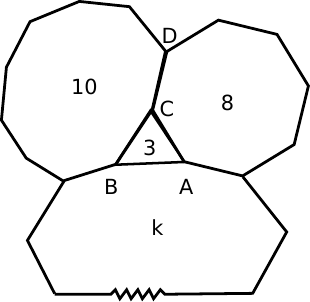}
(c)
\end{minipage}
\begin{minipage}[b]{0.4\textwidth}
\centering
\includegraphics[width=4.2cm,height=4.2cm]{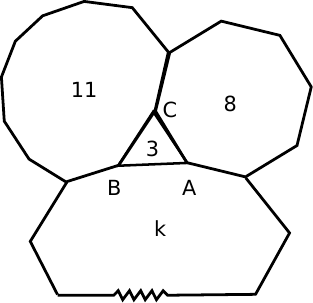}
(d)
\end{minipage}
\begin{minipage}[b]{0.4\textwidth}
\centering
\includegraphics[width=4.2cm,height=4.2cm]{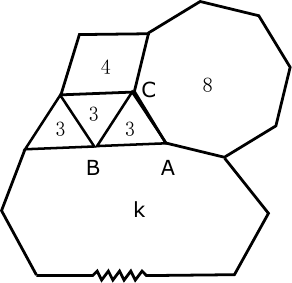}
(e)
\end{minipage}
\begin{minipage}[b]{0.4\textwidth}
\centering
\includegraphics[width=4.2cm,height=4.2cm]{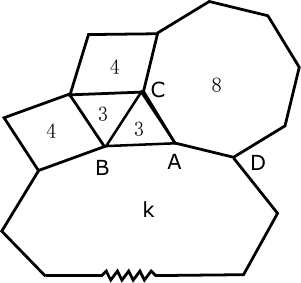}
(f)
\end{minipage}
\begin{minipage}[b]{0.4\textwidth}
\centering
\includegraphics[width=4.2cm,height=4.2cm]{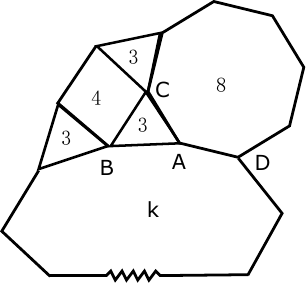}
(g)
\end{minipage}
\begin{minipage}[b]{0.4\textwidth}
\centering
\includegraphics[width=4.2cm,height=4.2cm]{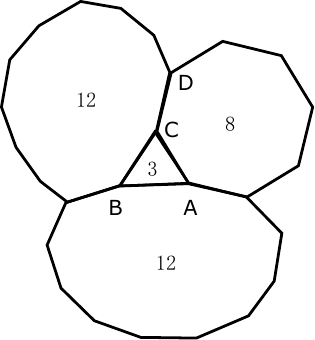}
(h)
\end{minipage}
\caption{\small Case 6.}
 \label{case6}
\end{figure}
 
\begin{itemize}
\item 
If $B=(3,8,k)$ for $k\leq 23,$ see Figure \ref{case6}(a), then $\Phi(C)=\frac{1}{12}$ and $\Phi(A)+\Phi(C)>\frac{1}{12}$.

\item
If $B=(3,9,k)$ for $k\leq  18 , $  see Figure \ref{case6}(b), then $\Phi(B)=\frac{1}{k}-\frac{1}{18}$ and $\Phi(C)=\frac{5}{72}$. Note that $D=(3,8,9)$ and $\Phi(D)=\frac{5}{72}$, which yields that $\Phi(A)+\Phi(B)+\Phi(C)+\Phi(D)>\frac{1}{12}$.

\item
If $B=(3,10,k)$ for $k\leq  15,$ see Figure \ref{case6}(c),  then $\Phi(B)=\frac{1}{k}-\frac{1}{15}$ and $\Phi(C)=\frac{7}{120}$.
For $D=(3,8,10),$ $\Phi(D)=\frac{7}{120}$ and $\Phi(A)+\Phi(B)+\Phi(C)+\Phi(D)>\frac{1}{12}$.
 
 \item
If $B=(3,11,k)$ for $k\leq  13 ,$ see Figure \ref{case6}(d),  then $\Phi(B)=\frac{1}{k}-\frac{5}{66}$ and $\Phi(C)=\frac{13}{264}$.
Hence $\Phi(A)+\Phi(B)+\Phi(C)>\frac{1}{12}$.
 
 \item
If $B=(3,3,3,k), $ see Figure \ref{case6}(e), then $\Phi(B)=\frac{1}{k}$.  The only nontrivial  case for $C$ is given by
$C=(3,3,4,8)$, which yields $\Phi(C)=\frac{1}{24}$. 
Since $k\leq 23$, $\Phi(A)+\Phi(B)+\Phi(C)>\frac{1}{12}$. 
 
 \item
If $B=(3,3,4,k)$ for $k\leq12$, then $\Phi(B)=\frac{1}{k}-\frac{1}{12}$. There are two subcases: The first one is shown in Figure \ref{case6}(f).
The only nontrivial case is
$C=(3,3,4,8)$, which implies $\Phi(C)=\frac{1}{24}$. If $k<12$, then $\Phi(A)+\Phi(B)+\Phi(C)>\frac{1}{12}$.
If $k=12$, note that $\Phi(D)=\frac{1}{24}$, which yields $\Phi(A)+\Phi(C)+\Phi(D)>\frac{1}{12}$;
The second case is depicted in Figure \ref{case6}(g). Similarly, we also have the total curvature $\Phi(G)>\frac{1}{12}$.

 \item
If $B=(3,12,12),  $ see Figure \ref{case6}(h), then  $\Phi(A)=\frac{1}{24}$, $\Phi(B)=0$, and $\Phi(C)=\frac{1}{24}$. Since $\Phi(D)=\frac{1}{24}$, $\Phi(A)+\Phi(C)+\Phi(D)>\frac{1}{12}$.  

\end{itemize}  

\item[Case 7]
 $A=(3,9,k), 9\leq k\leq17.$ In this case, $\Phi(A)=\frac{1}{k}-\frac{1}{18}$.    The possible patterns of  $B$ are
   $ (3,9,k), k\leq 17, $  $(3,10,k),k\leq 15,$
 $(3,11,k),k\leq 13, $ $(3,3,3,k),$ $(3,3,4,k), k\leq12$, $(3,9,18)$ and $(3,12,12)$.
 
 \begin{figure}[tb]
\begin{minipage}[b]{0.4\textwidth}
\centering
\includegraphics[width=4.2cm,height=4.2cm]{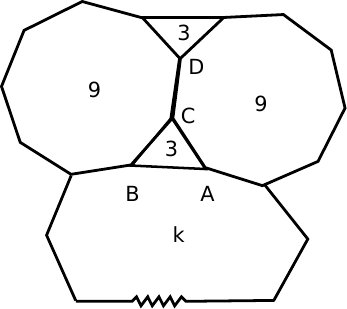}
(a)
\end{minipage}
\begin{minipage}[b]{0.4\textwidth}
\centering
\includegraphics[width=4.2cm,height=4.2cm]{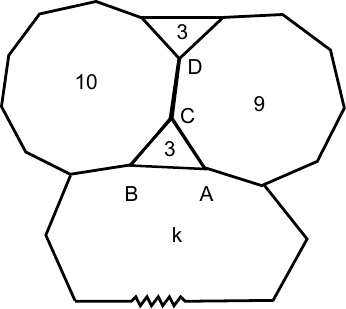}
(b)
\end{minipage}
\begin{minipage}[b]{0.4\textwidth}
\centering
\includegraphics[width=4.2cm,height=4.2cm]{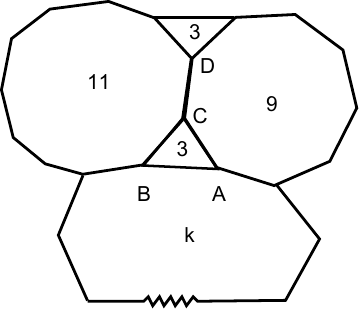}
(c)
\end{minipage}
\begin{minipage}[b]{0.4\textwidth}
\centering
\includegraphics[width=4.2cm,height=4.2cm]{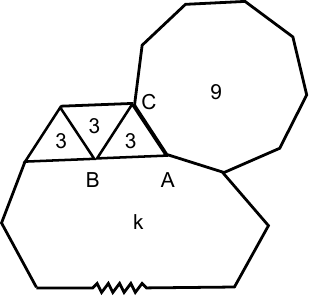}
(d)
\end{minipage}
\begin{minipage}[b]{0.4\textwidth}
\centering
\includegraphics[width=4.2cm,height=4.2cm]{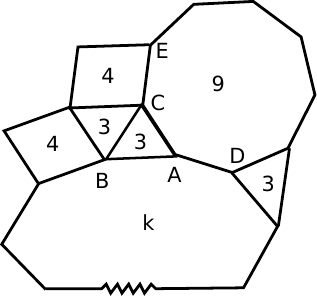}
(e)
\end{minipage}
\begin{minipage}[b]{0.4\textwidth}
\centering
\includegraphics[width=4.2cm,height=4.2cm]{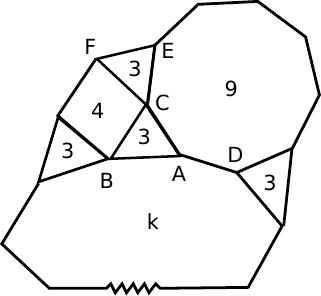}
(f)
\end{minipage}
\begin{minipage}[b]{0.4\textwidth}
\centering
\includegraphics[width=4.2cm,height=4.2cm]{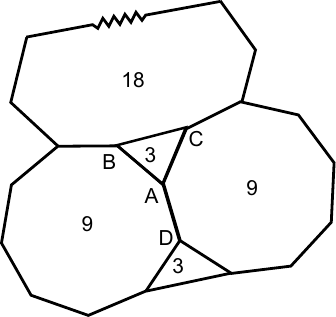}
(g)
\end{minipage}
\begin{minipage}[b]{0.4\textwidth}
\centering
\includegraphics[width=4.2cm,height=4.2cm]{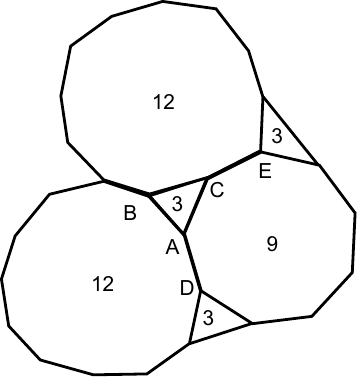}
(i)
\end{minipage}
\caption{\small Case 7.}
 \label{case7}
\end{figure}
 
\begin{itemize}
\item 
If $B=(3,9,k)$ for $k\leq17,$ see Figure \ref{case7}(a), then  $\Phi(C)=\frac{1}{18}$ and $\Phi(D)=\frac{1}{18}$.  This implies $\Phi(C)+\Phi(D)=\frac{2}{k}>\frac{1}{12}$. Hence $\Phi(A)+\Phi(C)+\Phi(D)>\frac{1}{12}$.
 
\item 
If $B=(3,10,k)$ for $k\leq15,$  see Figure \ref{case7}(b), then $\Phi(C)=\frac{2}{45}$ and $\Phi(D)=\frac{2}{45}$. This implies $\Phi(C)+\Phi(D)>\frac{1}{12}$. Hence $\Phi(A)+\Phi(C)+\Phi(D)>\frac{1}{12}$.
 
 \item 
If $B=(3,11,k)$ for $k\leq 13, $  see Figure \ref{case7}(c), then $\Phi(B)=\frac{1}{k}-\frac{5}{66}$, $\Phi(C)=\frac{7}{198}$, and $\Phi(D)=\frac{7}{198}$. Hence $\Phi(A)+\Phi(B)+\Phi(C)+\Phi(D)> \frac{1}{12}$.
 
 \item 
If $B=(3,3,3,k)$ for $k\leq 17,$ see Figure \ref{case7}(d), then $\Phi(B)=\frac{1}{k}$ and $\Phi(C)\geq \frac{1}{36}$.
Hence $\Phi(A)+\Phi(B)+\Phi(C)> \frac{1}{12}$.
 
 \item 
If $B=(3,3,4,k)$ for $k\leq12,$ then we have two subcases: The first one is shown in Figure \ref{case7}(e).  We have $\Phi(B)=\frac{1}{k}-\frac{1}{12}$, $\Phi(C)\geq \frac{1}{36}$, and $\Phi(D)=\frac{1}{k}-\frac{1}{18}$. By $\Phi(E)\geq\frac{1}{252}$, $\Phi(A)+\Phi(B)+\Phi(C)+\Phi(D)+\Phi(E)>\frac{1}{12}$;
The second one is shown in Figure \ref{case7}(f), where $\Phi(B)=\frac{1}{k}-\frac{1}{12}$, $\Phi(C)=\frac{1}{36}$, and
$\Phi(D)=\frac{1}{k}-\frac{1}{18}$. Note that if $\Phi(E)\neq0$, then $\Phi(E)\geq\frac{1}{306}$ and if $\Phi(E)=0$, then $E=(3,9,18)$, which yields $F=(3,4,18)$ and $\Phi(F)=\frac{5}{36}$. So we always have $\Phi(E)+\Phi(F)\geq\frac{1}{306}$, hence $\Phi(A)+\Phi(B)+\Phi(C)+\Phi(D)+\Phi(E)+\Phi(F)>\frac{1}{12}$. 
 
 \item 
If $B=(3,9,18),$ see Figure \ref{case7}(g), then $\Phi(A)=\frac{1}{18}$, $\Phi(B)=0$, $\Phi(C)=0$,
  and $\Phi(D)=\frac{1}{18}$. Hence $\Phi(A)+\Phi(D)>\frac{1}{12}$.
 
 
 \item 
If $B=(3,12,12), $  see Figure \ref{case7}(i), then  $\Phi(A)=\frac{1}{36}$, $\Phi(B)=0$, and $\Phi(C)=\Phi(D)=\Phi(E)=\frac{1}{36}.$ This yields that $\Phi(A)+\Phi(C)+\Phi(D)+\Phi(E)>\frac{1}{12}$.
 
%

 \end{itemize}

\item[Case 8]
 $A=(3,10,k)$ for $10\leq k \leq 14.$ In this case, $\Phi(A)=\frac{1}{k}-\frac{1}{15}$. The possible patterns of  $B$ are
   $(3,10,k),k\leq 14,$
 $(3,11,k),k\leq 13, $ $(3,3,3,k),$ $(3,3,4,k), k\leq 11,$  $(3,10,15),$ $(3,12,12) $ and $(3,3,4,12)$.

 \begin{figure}[tb]
\begin{minipage}[b]{0.4\textwidth}
\centering
\includegraphics[width=4.2cm,height=4.2cm]{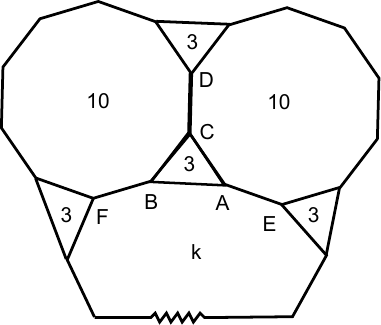}
(a)
\end{minipage}
\begin{minipage}[b]{0.4\textwidth}
\centering
\includegraphics[width=4.2cm,height=4.2cm]{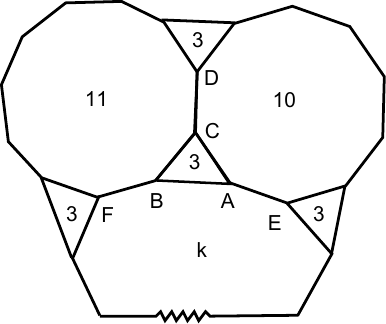}
(b)
\end{minipage}
\begin{minipage}[b]{0.4\textwidth}
\centering
\includegraphics[width=4.2cm,height=4.2cm]{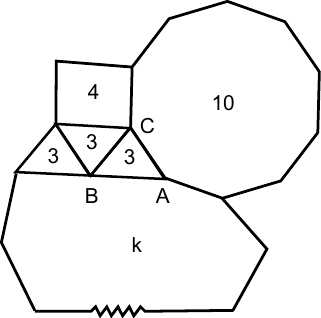}
(c)
\end{minipage}
\begin{minipage}[b]{0.4\textwidth}
\centering
\includegraphics[width=4.2cm,height=4.2cm]{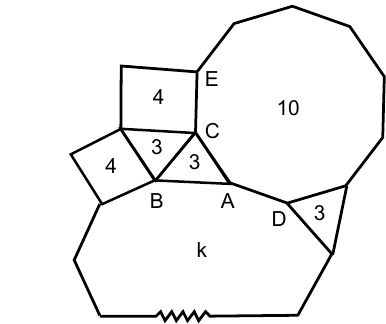}
(d)
\end{minipage}
\begin{minipage}[b]{0.4\textwidth}
\centering
\includegraphics[width=4.2cm,height=4.2cm]{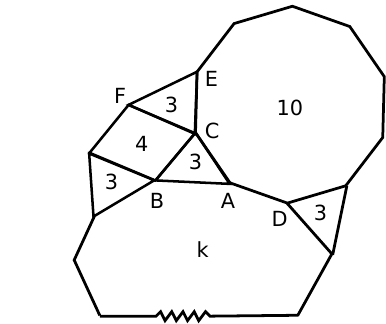}
(e)
\end{minipage}
\caption{\small Case 8: (a)-(e).}
 \label{case8: (a)-(e).}
\end{figure}

 \begin{figure}[tb]
\begin{minipage}[b]{0.4\textwidth}
\centering
\includegraphics[width=4.2cm,height=4.2cm]{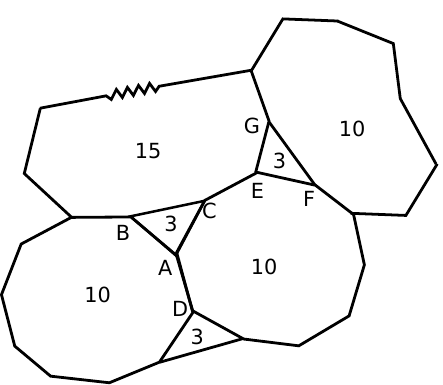}
(f)
\end{minipage}
\begin{minipage}[b]{0.4\textwidth}
\centering
\includegraphics[width=4.2cm,height=4.2cm]{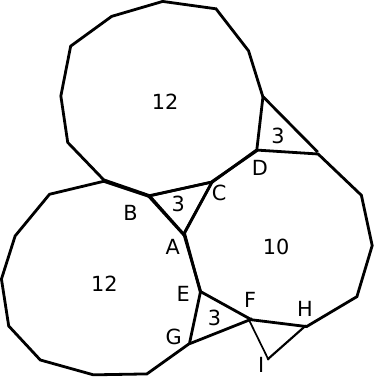}
(g)
\end{minipage}
\begin{minipage}[b]{0.4\textwidth}
\centering
\includegraphics[width=4.2cm,height=4.2cm]{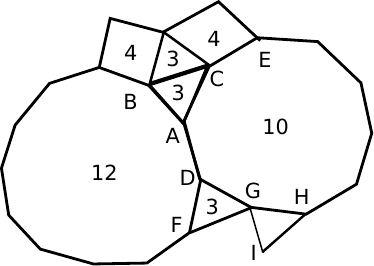}
(h)
\end{minipage}
\begin{minipage}[b]{0.4\textwidth}
\centering
\includegraphics[width=4.2cm,height=4.2cm]{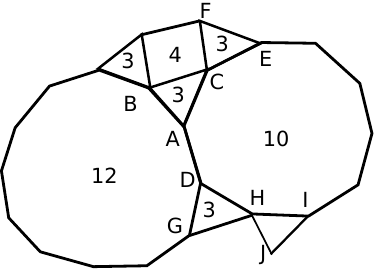}
(i)
\end{minipage}
\caption{\small Case 8 (f)-(i).}
 \label{case8: (f)-(i).}
\end{figure}

\begin{itemize}
\item 
If $B=(3,10,k)$ for $k\leq 14,$ see Figure \ref{case8: (a)-(e).}(a), then $\Phi(B)=\frac{1}{k}-\frac{1}{15}$, $\Phi(C)=\frac{1}{30}$, $\Phi(D)=\frac{1}{30}$, $\Phi(E)=\frac{1}{k}-\frac{1}{15},$ and
 $\Phi(F)=\frac{1}{k}-\frac{1}{15}$. Since $k\leq 14$,  
 $\Phi(A)+\Phi(B)+\Phi(C)+\Phi(D)+\Phi(E)+\Phi(F)>\frac{1}{12}$.
 
\item 
If $B=(3,11,k)$ for $k\leq 13,$ see Figure \ref{case8: (a)-(e).}(b), then $\Phi(B)=\frac{1}{k}-\frac{5}{66}$, $\Phi(C)=\frac{4}{165}$,
  $\Phi(D)=\frac{4}{165}$, $\Phi(E)=\frac{1}{k}-\frac{1}{15},$ and
 $\Phi(F)=\frac{1}{k}-\frac{5}{66}$. If $ k\leq 12$, then
 $\Phi(A)+\Phi(B)+\Phi(C)+\Phi(D)+\Phi(E)+\Phi(F)>\frac{1}{12}$. For $k=13$, since the curvature of any other vertex on the tridecagon and hendecagon is at least $\frac{1}{858}$ and there are $16$ such vertices, $\Phi(G)\geq\Phi(A)+\Phi(B)+\Phi(C)+\Phi(D)+\Phi(E)+\Phi(F)+\frac{16}{858}>\frac{1}{12}$.

 \item 
 If $B=(3,3,3,k)$ for $k\leq 14,$ see Figure \ref{case8: (a)-(e).}(c), then $\Phi(B)=\frac{1}{k}$ and $\Phi(C)\geq \frac{1}{60}$. Hence $\Phi(A)+\Phi(B)+\Phi(C)>\frac{1}{12}$.
 
 \item 
If $B=(3,3,4,k), k\leq 11,$  then we have two subcases: The first one is shown in Figure \ref{case8: (a)-(e).}(d). We have $\Phi(B)=\frac{1}{k}-\frac{1}{12}$, $\Phi(C)\geq \frac{1}{60}$,
$\Phi(D)=\frac{1}{k}-\frac{1}{15}$, and $\Phi(E)\geq \frac{1}{60}.$ Hence $\Phi(A)+\Phi(B)+\Phi(C)+\Phi(D)+\Phi(E)>\frac{1}{12}$. The other is depicted in Figure \ref{case8: (a)-(e).}(e), in which $\Phi(B)=\frac{1}{k}-\frac{1}{12}$, $\Phi(C)=\frac{1}{60}$ and
$\Phi(D)=\frac{1}{k}-\frac{1}{15}$. Now, we consider the vertex $E$. The possible patterns of $E$ are $(3,10,l),10\leq l\leq15$, $(3,3,3,10)$ and $(3,3,4,10)$. If $E=(3,10,l), l\leq12, (3,3,3,10)$ or $(3,3,4,10)$, then $\Phi(E)\geq\frac{1}{60}$ and $\Phi(A)+\Phi(B)+\Phi(C)+\Phi(D)+\Phi(E)>\frac{1}{12}$. If $E=(3,10,l),12<l\leq15$, then $F=(3,4,l)$ and $\Phi(F)>\frac{1}{12}$. Hence $\Phi(A)+\Phi(C)+\Phi(F)>\frac{1}{12}$.
 
 \item 
If $B=(3,10,15),$ see Figure \ref{case8: (f)-(i).}(f), then  $\Phi(A)=\frac{1}{30}$, $\Phi(B)=\Phi(C)=0$,
 $\Phi(D)=\frac{1}{30}$,  $E=(3,10,15)$ and $\Phi(E)=0$. We consider the vertex $G$.  The only nontrivial patterns for $G$ are $(3,10,15)$ and $(3,3,3,15)$. If $G=(3,10,15)$, then $F=(3,10,10)$ and $\Phi(F)=\frac{1}{30}$. Hence $\Phi(A)+\Phi(D)+\Phi(F)>\frac{1}{12}$. If $G=(3,3,3,15)$, then $\Phi(G)=\frac{1}{15}$ which yields $\Phi(A)+\Phi(D)+\Phi(G)>\frac{1}{12}$.

 \item 
If $B=(3,12,12), $ see Figure \ref{case8: (f)-(i).}(g), then  $\Phi(A)=\frac{1}{60}$, $\Phi(B)=0$, and $\Phi(C)=\Phi(D)=\Phi(E)=\frac{1}{60}$. The possible patterns of $G$ are
 $(3,10,12), $ $(3,11,12)$, $(3,3,3,12)$, $(3,12,12)$ and $(3,3,4,12)$. If $G=(3,10,12)$, then $\Phi(G)=\frac{1}{60}$, $F=(3,10,10)$ and $\Phi(F)=\frac{1}{30}$, which yields that $\Phi(A)+\Phi(C)+\Phi(D)+\Phi(E)+\Phi(F)>\frac{1}{12}$. If $G=(3,11,12)$, then $F=(3,10,11)$ and $\Phi(F)=\frac{4}{165}$. Hence $\Phi(A)+\Phi(C)+\Phi(D)+\Phi(E)+\Phi(F)>\frac{1}{12}$. If $G=(3,3,3,12)$, then $\Phi(G)=\frac{1}{12}$, which implies $\Phi(A)+\Phi(E)+\Phi(G)>\frac{1}{12}$. If $G=(3,12,12)$, then $F=(3,10,12)$ and $\Phi(F)=\Phi(H)=\frac{1}{60}$. This yields $\Phi(A)+\Phi(C)+\Phi(D)+\Phi(E)+\Phi(F)+\Phi(H)>\frac{1}{12}$. If $G=(3,3,4,12)$, then the pattern of $F$ is either $(3,3,3,10)$ or $(3,3,4,10).$ The former case is trivial since $\Phi(F)>\frac{1}{12}.$ So that $F=(3,3,4,10)$ and $\Phi(F)=\frac{1}{60}$. Moreover, if $\Phi(H)=0,$ that is $H=(3,10,15)$, then we have $I=(3,4,15)$ and $\Phi(I)>\frac{1}{12}$. If $\Phi(H)\neq0$, then $\Phi(H)\geq\frac{1}{210}$, which yields $\Phi(A)+\Phi(C)+\Phi(D)+\Phi(E)+\Phi(F)+\Phi(H)>\frac{1}{12}$.
{Note that in this subcase, we have proven that 
\begin{equation}\label{eq:fghi}
\Phi(F)+\Phi(G)+\Phi(H)+\Phi(I)\geq\frac{3}{140}.
\end{equation}
}
 
  \item 
If $B=(3,3,4,12),$ then we have two subcases: The first one is shown in Figure \ref{case8: (f)-(i).}(h). We have $\Phi(A)=\frac{1}{60}$, $\Phi(B)=0$, $\Phi(C)\geq\frac{1}{60}$,
$\Phi(D)=\frac{1}{60}$ and $\Phi(E)\geq \frac{1}{60}$. {Similar to \eqref{eq:fghi},} we have $\Phi(F)+\Phi(G)+\Phi(H)+\Phi(I)\geq\frac{3}{140}$.
Hence $\Phi(A)+\Phi(B)+\Phi(C)+\Phi(D)+\Phi(E)+\Phi(F)+\Phi(G)+\Phi(H)+\Phi(I)>\frac{1}{12}$;
The other is depicted in Figure \ref{case8: (f)-(i).}(i) where  $\Phi(A)=\frac{1}{60}$, $\Phi(B)=0$ and $\Phi(C)=\Phi(D)=\frac{1}{60}$. {Similar to \eqref{eq:fghi},} $\Phi(G)+\Phi(H)+\Phi(I)+\Phi(J)\geq\frac{3}{140}$. 
Now, we consider the vertex $E$. The possible patterns of $E$ are $(3,10,l),10\leq l\leq15$, $(3,3,3,10)$ and $(3,3,4,10)$. If $E=(3,10,l),l\leq12, (3,3,3,10)$ or $(3,3,4,10)$, then $\Phi(E)\geq\frac{1}{60}$. Hence $\Phi(A)+\Phi(C)+\Phi(E)+\Phi(D)+\Phi(G)+\Phi(H)+\Phi(I)+\Phi(J)>\frac{1}{12}$. If $E=(3,10,l),12<l\leq15$, then $F=(3,4,l)$, which yields $\Phi(F)>\frac{1}{12}$. Hence $\Phi(A)+\Phi(C)+\Phi(F)>\frac{1}{12}$.

 \end{itemize}

\item[Case 9]
$A=(3,11,k)$ for $11\leq k \leq 13.$ In this case, $\Phi(A)=\frac{1}{k}-\frac{5}{66}$. The possible patterns of  $B$ are $(3,11,k), $ $(3,3,3,k),$  $(3,3,4,11),$ $(3,12,12)$ and $(3,3,4,12)$.

 \begin{figure}[tb]
\begin{minipage}[b]{0.4\textwidth}
\centering
\includegraphics[width=4.2cm,height=4.2cm]{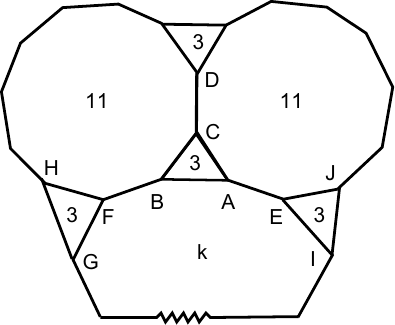}
(a)
\end{minipage}
\begin{minipage}[b]{0.4\textwidth}
\centering
\includegraphics[width=4.2cm,height=4.2cm]{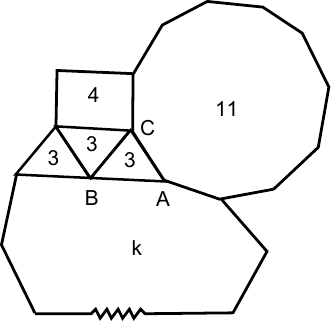}
(b)
\end{minipage}
\begin{minipage}[b]{0.4\textwidth}
\centering
\includegraphics[width=4.2cm,height=4.2cm]{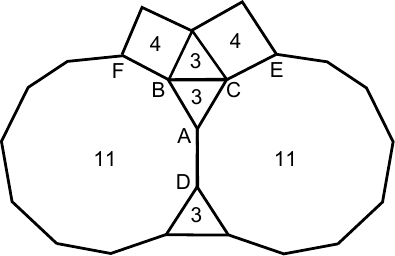}
(c)
\end{minipage}
\begin{minipage}[b]{0.4\textwidth}
\centering
\includegraphics[width=4.2cm,height=4.2cm]{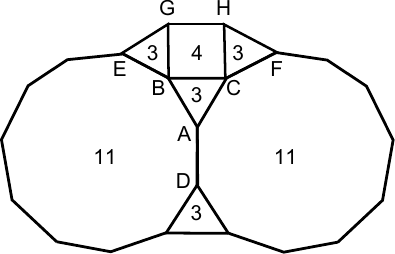}
(d)
\end{minipage}
\begin{minipage}[b]{0.4\textwidth}
\centering
\includegraphics[width=4.2cm,height=4.2cm]{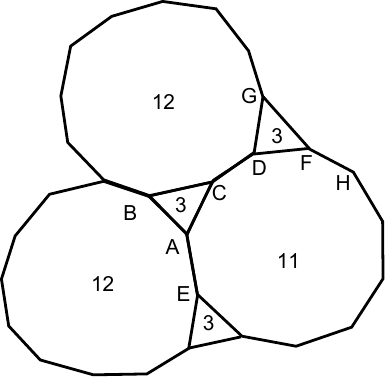}
(e)
\end{minipage}
\caption{\small Case 9.}
 \label{case9}
\end{figure}

\begin{itemize}
\item 
If $B=(3,11,k)$ for $11\leq k \leq 13,$  see Figure \ref{case9}(a),  then $\Phi(B)=\frac{1}{k}-\frac{5}{66}$, $\Phi(C)=\Phi(D)=\frac{1}{66}$, and $\Phi(E)=\Phi(F)=\frac{1}{k}-\frac{5}{66}$. 
For $k=11 $, $\Phi(A)+\Phi(B)+\Phi(C)+\Phi(D)+\Phi(E)+\Phi(F)> \frac{1}{12}$.
We consider the case of $k=12$. If none of the two hendecagons is adjacent to a tridecagon, then the curvature of each vertex on the two hendecagons is at least $\frac{1}{132}$, which implies that $\Phi(G)\geq(11+11-2)\times\frac{1}{132}=\frac{1}{12}+\frac{3}{44}>\frac{1}{12}$. If at least one of the hendecagons is adjacent a tridecagon, then both the vertices on the two hendecagons, except $A,B,C,D,E$ and $F$, and the vertices on the tridecagon are of the curvature at least $\frac{1}{858}$. Since there are at least $7+7+13-4=23$ such vertices, 
 $\Phi(A)+\Phi(B)+\Phi(C)+\Phi(D)+\Phi(E)+\Phi(F)+\frac{23}{858}>\frac{1}{12}$.
 
If $k=13$, we consider the vertex $G$. The nontrivial patterns of $G$ are $(3,11,13)$ and $(3,3,3,13)$. If $G=(3,11,13)$, then $H=(3,11,11)$. Hence, $\Phi(G)+\Phi(H)=\frac{1}{858}+\frac{1}{66}=\frac{7}{429}$. If $G=(3,3,3,13)$, then we get $\Phi(G)=\frac{1}{12}>\frac{7}{429}$. So in both cases, we always have $\Phi(G)+\Phi(H)\geq\frac{7}{429}$. Similarly, $\Phi(I)+\Phi(J)\geq \frac{7}{429}$.
And the curvatures of other vertices (except $A,B,\cdots,J$) of the two hendecagons and the tridecagon are at least $\frac{1}{858}$. Since there are $6+6+7=19$ such vertices,
$\Phi(A)+\Phi(B)+\Phi(C)+\Phi(D)+\Phi(E)+\Phi(F)+\Phi(G)+\Phi(H)+\Phi(I)+\Phi(J)+\frac{19}{858}>\frac{1}{12}$.
 
\item 
If $B=(3,3,3,k)$ for $11\leq k\leq 13,$ see Figure \ref{case9}(b), then $\Phi(B)=\frac{1}{k}$ and $\Phi(C)\geq \frac{1}{132}$.
Hence  $\Phi(A)+\Phi(B)+\Phi(C)>\frac{1}{12}$.
 
\item 
If $B=(3,3,4,11),$ then we have two subcases: The first one is shown in Figure \ref{case9}(c). We have $\Phi(A)=\frac{1}{66}$, $\Phi(B)=\frac{1}{132}$ and $\Phi(D)=\frac{1}{66}$. If $C=(3,3,3,11)$, then $\Phi(C)=\frac{1}{11}>\frac{1}{12}$, which implies $\Phi(A)+\Phi(C)>\frac{1}{12}$. So we only need to consider $C=(3,3,4,11).$ In that case, $\Phi(C)=\frac{1}{132},$ $\Phi(E)\geq\frac{1}{132}$ and $\Phi(F)\geq\frac{1}{132}$. If none of the two hendecagons is adjacent to a tridecagon, then the curvatures of the vertices on the two hendecagons,  except $A,B,C,D,E$ and $F$, are at least $\geq\frac{1}{132}$, which yields
$\Phi(G)\geq\Phi(A)+\Phi(B)+\Phi(C)+\Phi(D)+\Phi(E)+\Phi(F)+\frac{14}{132}>\frac{1}{12}$. If at least one of the hendecagons is adjacent to a tridecagon, then the curvatures of the vertices on the two hendecagons (except $A,B,C,D,E,F$) and the tridecagon are at least $\frac{1}{858}.$ Hence
$\Phi(G)\geq\Phi(A)+\Phi(B)+\Phi(C)+\Phi(D)+\Phi(E)+\Phi(F)+\frac{23}{858}>\frac{1}{12}$. This proves the first case; The second one is shown in Figure \ref{case9}(d). We have $\Phi(A)=\frac{1}{66}$, $\Phi(B)=\frac{1}{132}$, $\Phi(C)=\frac{1}{132}$ and $\Phi(D)=\frac{1}{66}$. If $E=(3,11,13)$, then $G=(3,4,13)$, which implies $\Phi(G)=\frac{25}{156}>\frac{1}{12}.$ Hence $\Phi(A)+\Phi(B)+\Phi(G)>\frac{1}{12}$. If $E\neq(3,11,13)$, then $\Phi(E)\geq\frac{1}{132}$. Similarly we have $\Phi(F)\geq\frac{1}{132}$. If none of the two hendecagons is adjacent to a tridecagon, then the curvatures of the vertices on the two hendecagons, except $A,B,C,D,E$ and $F$, are at least $\frac{1}{132}$, which yields
$\Phi(G)\geq\Phi(A)+\Phi(B)+\Phi(C)+\Phi(D)+\Phi(E)+\Phi(F)+\frac{14}{132}>\frac{1}{12}$. If at least one of the hendecagons is adjacent to a tridecagon, then the curvatures of the vertices on the two hendecagons, except $A,B,C,D,E,F$, and the vertices on the tridecagon are at least $\frac{1}{858}.$ Hence $\Phi(G)\geq\Phi(A)+\Phi(B)+\Phi(C)+\Phi(D)+\Phi(E)+\Phi(F)+\frac{23}{858}>\frac{1}{12}$.

\item 
If $B=(3,12,12),$ see Figure \ref{case9}(e), then we can construct a graph with $B=(3,12,12)$ and the total curvature $\frac{1}{12}$ {(see Figure \ref{fig5})}. We need to prove that in other cases the total curvature is no less than $\frac{1}{12}$. Note that $\Phi(A)=\frac{1}{132}$, $\Phi(B)=0$ and $\Phi(C)=\Phi(D)=\Phi(E)=\frac{1}{132}$. Since the pattern of $F$ cannot be $(3,11,13)$ (otherwise, $G$ would be $(3,12,13)$ whose curvature is negative), $\Phi(F)\geq\frac{1}{132}$. So the possible patterns of $F$ are $(3,11,l),11\leq l\leq12$, $(3,3,3,11)$ and $(3,3,4,11)$. For any case of $F$, $H$ cannot be adjacent to a tridecagon (otherwise, either there exists a vertex with negative curvature or the sum of curvature is greater than $\frac{1}{12}$). So $\Phi(H)\geq\frac{1}{132}$. It works similarly for other vertices on the hendecagon, and we get the curvature of each vertex of the hendecagon is at least $\frac{1}{132}$. Hence the total curvature is at least $\frac{1}{12}$ and the total curvature is $\frac{1}{12}$ if and only if the curvature of each vertex of the hendecagon is $\frac{1}{132}$ and other vertices have vanishing curvature.
  
 \item 
If $B=(3,3, 4,12), $ then we can construct a graph with $B=(3,3,4,12)$ and total curvature $\frac{1}{12}$ {(see also Figure \ref{fig5})}. We only need to show that in other cases the total curvature is at least $\frac{1}{12}$. The proof is similar to the above case with $B=(3,12,12)$ and is omitted here.

 \end{itemize}
\item[Case 10]
$A=(4,4,k).$ In this case, $\Phi(A)=\frac{1}{k}$. For $k=12$, we can construct a graph with total curvature $\frac{1}{12}$ {see Figure \ref{4-4-12}}. So we only need to consider $k>12$. Denote by $B$ and $C$ the neighbors of $A$ incident to the $k$-gon. In this case, the possible patterns of $B$ and $C$ are 
$(4,4,k)$ and $(4,5,k), k\leq20$.

 \begin{figure}[tb]
\begin{minipage}[b]{0.4\textwidth}
\centering
\includegraphics[width=4.2cm,height=4.2cm]{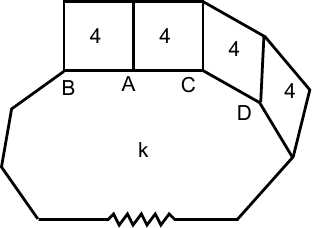}
(a)
\end{minipage}
\begin{minipage}[b]{0.4\textwidth}
\centering
\includegraphics[width=4.2cm,height=4.2cm]{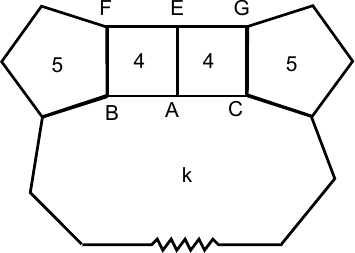}
(b)
\end{minipage}
\caption{\small Case 10.}
 \label{case10}
\end{figure}

\begin{itemize}
\item   
If $B=(4,4,k), $ see Figure \ref{case10}(a), then  $\Phi(B)=\frac{1}{k}$. 
Hence if $k\leq 23$, then $\Phi(A)+\Phi(B)>\frac{1}{12}$.  
 If  $k\geq 24$, the pattern of $C$ has to be $(4,4,k)$. Hence $\Phi(C)=\frac{1}{k}.$ 
Similarly $D=(4,4,k)$ and $\Phi(D)=\frac{1}{k}$.  
 So that $\Phi(A)+\Phi(B)+\Phi(C)+\Phi(D)=\frac{4}{k}>\frac{1}{12}$ by $k\leq42$.
 
 \item   
If $B=(4,5,k)$ for $k\leq20,$ see Figure \ref{case10}(b), then $\Phi(B)=\frac{1}{k}-\frac{1}{20}$. 
If $k\leq 14$, then $\Phi(A)+\Phi(B)>\frac{1}{12}$.  
So we consider $k\geq 15$. If $C=(4,4,k)$, then noting that $k\leq20$, we have $\Phi(A)+\Phi(B)+\Phi(C)>\frac{1}{12}$. If $C=(4,5,k)$, then $\Phi(C)=\frac{1}{k}-\frac{1}{20}.$
We consider the vertex $E.$ The possible patterns of $E$ are $(4,4,l), (l\geq 12),$ $(3,3,4,4)$,  $(3,4,4,4)$, $(3,4,4,5)$,
 $(3,4,4,6)$, $(4,4,4,4)$ and $(3,3,3,4,4)$.
If $E=(4,4,l)$, $(3,3,4,4)$ or $(3,4,4,4)$, then by $l\leq20,$ $\Phi(E)\geq \frac{1}{20},$ which implies
 $\Phi(A)+\Phi(B)+\Phi(C)+\Phi(E)>\frac{1}{12}.$
If $E=(3,4,4,5)$, then $\Phi(E)=\frac{1}{30}$ and $\Phi(F)$ (or $\Phi(G))=\frac{3}{20}$, which yields $\Phi(A)+\Phi(B)+\Phi(C)+\Phi(E) +\Phi(F)+\Phi(G)>\frac{1}{12}$.
If $E=(3,4,4, 6), (3,3,3,4,4), (4,4,4,4)$, then $ \Phi(E)=0$ and $\Phi(F)+\Phi(G)\geq\frac{1}{15}$. Hence
 $\Phi(A)+\Phi(B)+\Phi(C)+\Phi(E)+\Phi(F)+\Phi(G)>\frac{1}{12}.$
%

\end{itemize}
\item[Case 11]
$A=(4,5,k)$ for $5\leq k \leq 19.$ In this case, $\Phi(A)=\frac{1}{k}-\frac{1}{20}$. Since $\Phi(A)>\frac{1}{12}$ for $k\leq7$. We only need to consider $k\geq8$. We denote by $B$ the neighbor of $A$ which is incident to the square and the $k$-gon, from this case to Case 13. The possible patterns of $B$ are 
$(4,5,k), k\leq19$ $(4,6,k),k\leq11$ $(4,7,k),k\leq9$, $(3,3,4,k),k\leq11$ $(4,6,12),$ $(4,8,8)$ and $(3,3,4,12)$.

 \begin{figure}[tb]
\begin{minipage}[b]{0.4\textwidth}
\centering
\includegraphics[width=4.2cm,height=4.2cm]{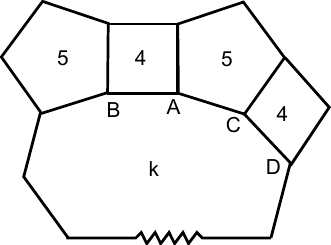}
(a)
\end{minipage}
\begin{minipage}[b]{0.4\textwidth}
\centering
\includegraphics[width=4.2cm,height=4.2cm]{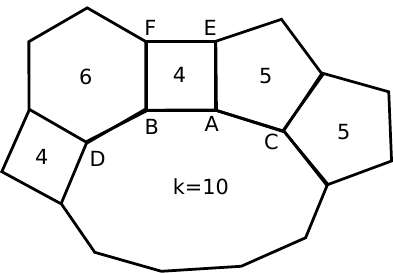}
(b)
\end{minipage}
\begin{minipage}[b]{0.4\textwidth}
\centering
\includegraphics[width=4.2cm,height=4.2cm]{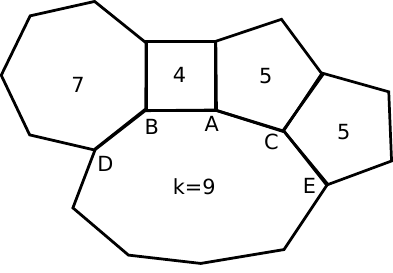}
(c)
\end{minipage}
\begin{minipage}[b]{0.4\textwidth}
\centering
\includegraphics[width=4.2cm,height=4.2cm]{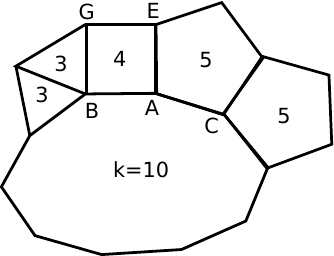}
(d)
\end{minipage}
\begin{minipage}[b]{0.4\textwidth}
\centering
\includegraphics[width=4.2cm,height=4.2cm]{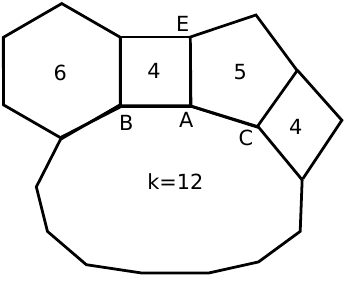}
(e)
\end{minipage}
\begin{minipage}[b]{0.4\textwidth}
\centering
\includegraphics[width=4.2cm,height=4.2cm]{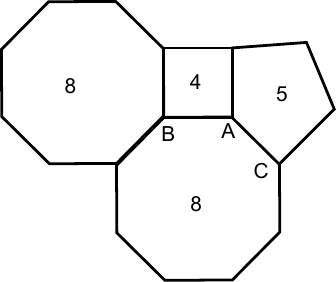}
(f)
\end{minipage}
\begin{minipage}[b]{0.4\textwidth}
\centering
\includegraphics[width=4.2cm,height=4.2cm]{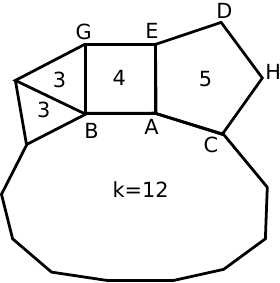}
(g)
\end{minipage}
\caption{\small Case 11.}
 \label{case11}
\end{figure}

\begin{itemize}
\item 
 If $B=(4,5,k)$, see Figure \ref{case11}(a), then $\Phi(B)=\frac{1}{k}-\frac{1}{20}$.  If $k\leq 10$, then $\Phi(A)+\Phi(B)>\frac{1}{12}$.
 It suffices to consider $k> 10.$ The possible pattern of $C$ is $(4,5,k)$ which implies $\Phi(C)=\frac{1}{k}-\frac{1}{20}$. If  $k\leq 12$, then $\Phi(A)+\Phi(B)+\Phi(C)>\frac{1}{12}.$ For $k>12$, the nontrivial pattern of $D$ is $(4,5,k)$ and $\Phi(D)=\frac{1}{k}-\frac{1}{20}$. We divide it into two subcases: If $k$ is even, then each vertex of the $k$-polygon, similar to the vertex $D$, has the curvature $\frac{1}{k}-\frac{1}{20}$; If $k$ is odd, then there exists a vertex on the $k$-polygon with a different pattern. However, the only nontrivial pattern except $(4,5,k)$ is $(3,3,3,k)$ (with the curvature $\frac{1}{k}$).
Hence for even $k$ ($k\leq18$), $\Phi(G)\geq k\times(\frac{1}{k}-\frac{1}{20})>\frac{1}{12}$; For odd $k$ ($k\leq19$), $\Phi(G)\geq(k-1)\times(\frac{1}{k}-\frac{1}{20})+\frac{1}{k}>\frac{1}{12}$.
%
%
%
%

\item 
 If $B=(4,6,k)$ for $6\leq k\leq 11,$ see Figure \ref{case11}(b), then $\Phi(B)=\frac{1}{k}-\frac{1}{12}$. 
 Note that if $k\leq 9$, then $\Phi(A)+\Phi(B)>\frac{1}{12}$. So we only need to consider $k\geq10$.
In this case, the possible pattern of $C$ is $(4,5,k)$ or $(5,5,10)$.
If $C=(4,5,k)$, then $\Phi(C)=\frac{1}{k}-\frac{1}{20}$. This yields $\Phi(A)+\Phi(B)+\Phi(C)>\frac{1}{12}$. 
If  $C=(5,5,10),$ then  $\Phi(A)=\frac{1}{20}$, $\Phi(B)=\frac{1}{60}$ and  $\Phi(C)=0$. The only nontrivial pattern for $D$ is $(4,6,10)$,
which implies $\Phi(D)=\frac{1}{60}$. Since $\Phi(E)\neq0$ (otherwise, $F=(4,6,20)$ whose curvature is negative), $\Phi(E)\geq\frac{1}{380}.$ Hence $\Phi(A)+\Phi(B)+\Phi(D)+\Phi(E)>\frac{1}{12}$.    
 
 \item  
 If $B=(4,7,k)$ for $k=8,9,$ see Figure \ref{case11}(c), then $\Phi(B)=\frac{1}{k}-\frac{3}{28}$. If $k=8$, then $\Phi(A)+\Phi(B)>\frac{1}{12}$.  
We only need to consider $k=9.$ The only nontrivial pattern of $D$ is $(4,7,9)$, which implies $\Phi(D)=\frac{1}{252}$. And the nontrivial patterns of $C$ are $(4,5,9)$ and $(5,5,9)$.
If $C=(4,5,9)$, then $\Phi(C)=\frac{11}{180}$. Hence $\Phi(A)+\Phi(B)+\Phi(C)+\Phi(D)>\frac{1}{12}$.
If $C=(5,5,9)$, then $\Phi(C)=\frac{1}{90}$ and $\Phi(E)\geq\frac{1}{90}$, which yields $\Phi(A)+\Phi(B)+\Phi(C)+\Phi(D)+\Phi(E)>\frac{1}{12}$.
 
 \item 
 If $B=(3,3,4,k)$ for $8\leq k\leq 11,$ see Figure \ref{case11}(d), then $\Phi(B)=\frac{1}{k}-\frac{1}{12}$.
 If $k\leq 9$, then $\Phi(A)+\Phi(B)>\frac{1}{12}$. So we only need to consider for $k= 10,11$. In these cases, the possible patterns of $C$ are $(4,5,k)$ and $(5,5,10)$.
 If $C=(4,5,k)$, then $\Phi(A)+\Phi(B)+\Phi(C)>\frac{1}{12}$.
 If $C=(5,5,10)$, then $\Phi(A)=\frac{1}{20},$  $\Phi(B)=\frac{1}{60},$ and $\Phi(C)=0$. Consider the vertex $E$. The possible patterns of $E$ are $(4,5,l),l\leq20, (3,3,4,5)$ and $(3,4,4,5)$.
If $E=(4,5,l)$ then $\Phi(E)=\frac{1}{l}-\frac{1}{20}$. For $l<15$, $\Phi(A)+\Phi(B)+\Phi(E)>\frac{1}{12}$. For $l\geq15$,  $G=(3,4,l)$ which implies $\Phi(G)>\frac{1}{12}.$ Hence $\Phi(A)+\Phi(B)+\Phi(G)>\frac{1}{12}$. 
If $E=(3,3,4,5)$ or $(3,4,4,5)$, then $\Phi(E)\geq\frac{1}{30}$. This yields $\Phi(A)+\Phi(B)+\Phi(E)>\frac{1}{12}$.
 
 \item 
If $B=(4,6,12),$ see Figure \ref{case11}(e), then  $\Phi(A)=\frac{1}{30},$ and $\Phi(B)=0$. The nontrivial pattern of $C$ is $(4,5,12),$ which implies $\Phi(C)=\frac{1}{30}$.
Consider the vertex $E$. The possible patterns of $E$ are $(4,5,l),l\leq12, (3,3,4,5)$ and $(3,4,4,5)$. In each case, we have $\Phi(E)\geq\frac{1}{30}$. Hence $\Phi(A)+\Phi(C)+\Phi(E)>\frac{1}{12}$.
 
 \item 
 If $B=(4,8,8), $ see Figure \ref{case11}(f), then $\Phi(A)=\frac{3}{40}$ and $\Phi(B)=0$. The nontrivial patterns of $C$ are $(4,5,8)$ and $(5,5,8)$, which implies $\Phi(C)\geq\frac{1}{40}.$ Hence $\Phi(A)+\Phi(C)>\frac{1}{12}$. 
 
\item 
If $B=(3,3,4,12)$, see Figure \ref{case11}(g), then  $\Phi(A)=\frac{1}{30}$ and $\Phi(B)=0$. The only nontrivial pattern of $C$ is $(4,5,12)$. Hence $\Phi(C)=\frac{1}{30}$. Consider the vertex $E$. The possible patterns of $E$ are $(4,5,l),l\leq20, (3,3,4,5)$ and $(3,4,4,5)$.
If $E=(4,5,l),$ then $\Phi(E)=\frac{1}{l}-\frac{1}{20}$. For $l<15$, $\Phi(A)+\Phi(C)+\Phi(E)>\frac{1}{12}$. For $15\leq l<20$,  $G=(3,4,l)$, which implies $\Phi(G)>\frac{1}{12}.$ Hence $\Phi(A)+\Phi(E)+\Phi(G)>\frac{1}{12}$. 
If $l=20$, then $\Phi(E)=0.$ The nontrivial pattern of $D$ is $(4,5,20)$. So that $H=(3,4,4,5)$ and $\Phi(H)=\frac{1}{30}$, which yields $\Phi(A)+\Phi(C)+\Phi(H)>\frac{1}{12}$.
If $E=(3,3,4,5)$ or $(3,4,4,5)$, then $\Phi(E)\geq\frac{1}{30}$. Hence $\Phi(A)+\Phi(B)+\Phi(E)>\frac{1}{12}$.

  \end{itemize}
\item[Case 12]
$A=(4,6,k)$ for $6\leq k \leq 11.$ In this case, $\Phi(A)=\frac{1}{k}-\frac{1}{12}$. Since for $k=6,11$ we can construct a graph with total curvature $\frac{1}{12}$ (see Figures \ref{fig3} and \ref{fig6}), it suffices to prove that for $A=(4,6,k),7\leq k\leq10$, $\Phi(G)>\frac{1}{12}$ and for $A=(4,6,11),$ $\Phi(G)\geq\frac{1}{12}$. The possible patterns of  $B$ are $(4,6,k), k\leq11,$ $(4,7,k), k\leq9$, $(3,3,4,k),k\leq11$ and $(4,8,8)$. 

 \begin{figure}[tb]
\begin{minipage}[b]{0.4\textwidth}
\centering
\includegraphics[width=4.2cm,height=4.2cm]{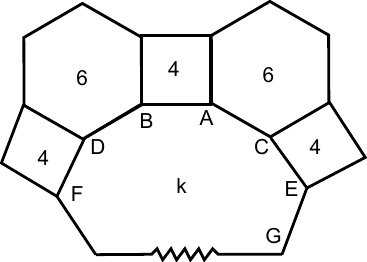}
(a)
\end{minipage}
\begin{minipage}[b]{0.4\textwidth}
\centering
\includegraphics[width=4.2cm,height=4.2cm]{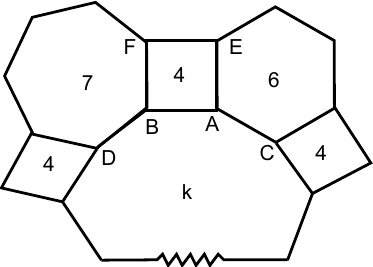}
(b)
\end{minipage}
\begin{minipage}[b]{0.4\textwidth}
\centering
\includegraphics[width=4.2cm,height=4.2cm]{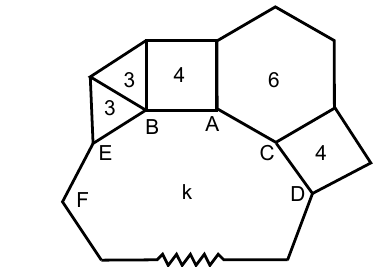}
(c)
\end{minipage}
\begin{minipage}[b]{0.4\textwidth}
\centering
\includegraphics[width=4.2cm,height=4.2cm]{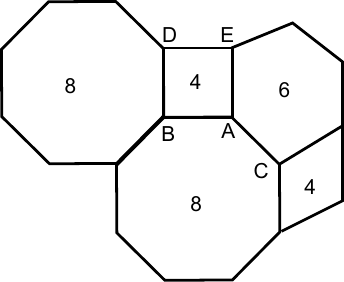}
(d)
\end{minipage}
\caption{\small Case 12.}
 \label{case12}
\end{figure}    
 
\begin{itemize}
\item 
If $B=(4,6,k)$ {for $7\leq k\leq 11$}, see Figure \ref{case12}(a), then $\Phi(B)=\frac{1}{k}-\frac{1}{12}$. For $k=7$, $\Phi(A)+\Phi(B)>\frac{1}{12}$. And for $k\geq8$, the only possible pattern of $C$ or $D$ is $(4,6,k)$. So that $\Phi(C)=\Phi(D)=\frac{1}{k}-\frac{1}{12}$. Hence, for $k\leq9$, $\Phi(A)+\Phi(B)+\Phi(C)+\Phi(D)>\frac{1}{12}$.
For $k=10,11$, we consider vertices $E,F$. The nontrivial possible patterns of $E,F$ are $(4,6,k)$ and $(3,3,4,k)$ which implies $\Phi(E)=\Phi(F)=\frac{1}{k}-\frac{1}{12}$. 
So for $k=10$, $\Phi(A)+\Phi(B)+\Phi(C)+\Phi(D)+\Phi(E)+\Phi(F)>\frac{1}{12}$.
For $k=11$, i.e. $E=(4,6,11)$ or $(3,3,4,11)$, we may argue as in Case 9 with $B=(3,12,12)$ to conclude $\Phi(G)\geq\frac{1}{12}$.
  
\item 
 If $B=(4,7,k)$ for $7\leq k\leq 9,$ see Figure \ref{case12}(b), then $\Phi(B)=\frac{1}{k}-\frac{3}{28}$. Since $\Phi(A)+\Phi(B)>\frac{1}{12}$ for $k=7,$
it suffices to consider the case $k=8,9$. Then the nontrivial pattern of $C$ is $(4,6,k),$ which implies $\Phi(C)=\frac{1}{k}-\frac{1}{12}$;
The nontrivial pattern of $D$ is $(4,7,k),$ which implies $\Phi(D)=\frac{1}{k}-\frac{3}{28}$. So that for $k=8$, $\Phi(A)+\Phi(B)+\Phi(C)+\Phi(D)>\frac{1}{12}$. For $k=9$, we consider the vertex $E$. The possible patterns of $E$ are $(4,6,l),l\leq12$, $(3,3,4,6)$ and $(3,4,4,6)$. 
If $E=(4,6,l)$, then $F=(4,7,l)$ ({$l\leq9$}). Hence $\Phi(E)=\frac{1}{l}-\frac{1}{12}\geq\frac{1}{36}$. This yields $\Phi(A)+\Phi(B)+\Phi(C)+\Phi(D)+\Phi(E)>\frac{1}{12}$.
If $E=(3,3,4,6)$, then $\Phi(A)+\Phi(E)>\frac{1}{12}$. For $E=(3,4,4,6)$, $F$ is either $(4,4,7)$ or $(3,3,4,7)$, which implies $\Phi(F)\geq\frac{5}{84}$. So that $\Phi(A)+\Phi(B)+\Phi(C)+\Phi(D)+\Phi(F)>\frac{1}{12}$.

\item 
If $B=(3,3,4,k)$ for $6\leq k\leq 11$, see Figure \ref{case12}(c), then $\Phi(B)=\frac{1}{k}-\frac{1}{12}$. Since $\Phi(A)+\Phi(B)>\frac{1}{12}$ for $k\leq 7,$ we consider $8\leq k\leq 11.$ For $k\geq8$, the pattern of $C$ is $(4,6,k)$, which implies $\Phi(C)=\frac{1}{k}-\frac{1}{12}$. In this case, $\Phi(A)+\Phi(B)+\Phi(C)>\frac{1}{12}$. For $k=9$, noting that $\Phi(D)\geq\frac{1}{252}$, we have $\Phi(A)+\Phi(B)+\Phi(C)+\Phi(D)>\frac{1}{12}$.
For $k=10$, the possible patterns of $D$ are $(4,6,k)$ and $(3,3,4,k)$, which implies $\Phi(D)=\frac{1}{k}-\frac{1}{12}$.
Now we consider the vertex $E.$ The possible patterns of $E$ are $(3,3,3,10)$, $(3,3,4,10)$ and $(3,10,15)$. If $E=(3,3,3,10),$ then $\Phi(A)+\Phi(B)+\Phi(E)>\frac{1}{12}$.  If $E=(3,3,4,10)$, then $\Phi(E)\geq\frac{1}{60}$ and $\Phi(F)\geq\frac{1}{306}$. Hence $\Phi(A)+\Phi(B)+\Phi(C)+\Phi(D)+\Phi(E)+\Phi(F)>\frac{1}{12}$.

If $k=11$, {we can construct a graph with total curvature $=\frac{1}{12}$, see Figure \ref{fig6}. } So, we only need to prove that in this case the total curvature is at least $\frac{1}{12}$. This can be obtained similarly as in the Case 9 with $B=(3,12,12)$.
 
\item 
 If $B=(4,8,8),$ see Figure \ref{case12}(d), then  $\Phi(A)=\frac{1}{24}$ and $\Phi(B)=0$. Since the nontrivial pattern of $C$ is $(4,6,8)$, which implies $\Phi(C)=\frac{1}{24}$ and $\Phi(A)+\Phi(C)=\frac{1}{12}$. We consider the vertex $E$. If $\Phi(E)\neq0$, then $\Phi(E)\geq\frac{1}{132}$, which yields $\Phi(A)+\Phi(C)+\Phi(E)>\frac{1}{12}$. If $\Phi(E)=0$, then $E=(3,4,4,6)$, so that $\Phi(D)\geq\frac{1}{24}$ and $\Phi(A)+\Phi(C)+\Phi(D)>\frac{1}{12}$.
 \end{itemize}

\item[Case 13]
$A=(4,7,k)$ for $7\leq k \leq 9.$ In this case, $\Phi(A)=\frac{1}{k}-\frac{3}{28}$. The possible patterns of $B$ are $(4,7,k)$, $(3,3,4,k)$ and $(4,8,8)$.

 \begin{figure}[tb]
\begin{minipage}[b]{0.4\textwidth}
\centering
\includegraphics[width=4.2cm,height=4.2cm]{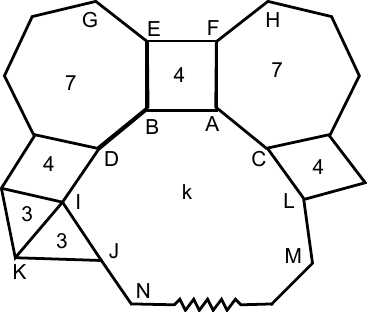}
(a)
\end{minipage}
\begin{minipage}[b]{0.4\textwidth}
\centering
\includegraphics[width=4.2cm,height=4.2cm]{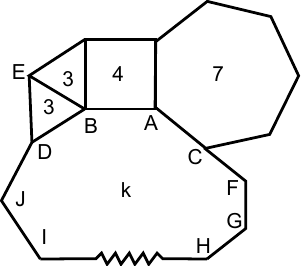}
(b)
\end{minipage}
\begin{minipage}[b]{0.4\textwidth}
\centering
\includegraphics[width=4.2cm,height=4.2cm]{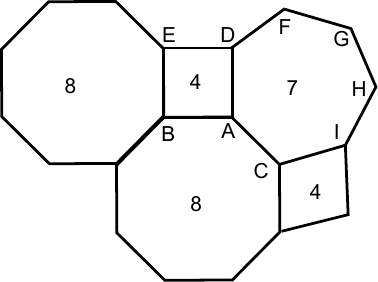}
(c)
\end{minipage}
\caption{\small Case 13.}
 \label{case13}
\end{figure}

\begin{itemize}
\item 
If $B=(4,7,k)$ for for $7\leq k \leq 9,$ see Figure \ref{case13}(a), then $\Phi(B)=\frac{1}{k}-\frac{3}{28}$.
The only nontrivial pattern of $C$ (or $D$) is $(4,7,k)$. So that $\Phi(C)=\frac{1}{k}-\frac{3}{28}$ and $\Phi(D)=\frac{1}{k}-\frac{3}{28}$.
Since $\Phi(A)+\Phi(B)+\Phi(C)+\Phi(D)>\frac{1}{12}$ for $k=7.$ We consider the case of $k=8$ or $9.$
Consider the vertex $E.$ The nontrivial patterns of $E$ are $(4,7,l)$ {$l\leq9$}, and $(3,3,4,7)$.
If $E=(3,3,4,7)$, then $F=(3,3,4,7)$, which yields $\Phi(A)+\Phi(E)+\Phi(F)>\frac{1}{12}$.
If $E=(4,7,l)$, then $F=(4,7,l),$ $\Phi(E)=\Phi(F)=\frac{1}{l}-\frac{3}{28}$. Furthermore, the pattern of $G$ (or $H$) has to be $(4,7,l)$, which implies $\Phi(G)=\Phi(H)=\frac{1}{l}-\frac{3}{28}$.
For $k=8$, by $l\leq9$ we have $\Phi(A)+\Phi(B)+\Phi(C)+\Phi(D)+\Phi(E)+\Phi(F)+\Phi(G)+\Phi(H)>\frac{1}{12}$.
For $k=9$, we consider the vertex $I$. The possible patterns of $I$ are $(4,7,9)$ and $(3,3,4,9)$. If $I=(3,3,4,9)$, then $\Phi(I)=\frac{1}{36}$. So that the possible patterns of $J$ are $(3,3,3,9)$, $(3,3,4,9)$ and $(3,9,18)$. If $J=(3,3,3,9)$ or $(3,3,4,9)$, then $\Phi(J)\geq\frac{1}{36}$, which yields 
$\Phi(A)+\Phi(B)+\Phi(C)+\Phi(D)+\Phi(E)+\Phi(F)+\Phi(G)+\Phi(H)+\Phi(I)+\Phi(J)>\frac{1}{12}$. And, if $J=(3,9,18)$, then $\Phi(K)\geq\frac{1}{18}$, which yields $\Phi(A)+\Phi(B)+\Phi(C)+\Phi(D)+\Phi(E)+\Phi(F)+\Phi(G)+\Phi(H)+\Phi(I)+\Phi(K)>\frac{1}{12}$.
If $I=(4,7,9)$, then $\Phi(I)=\frac{1}{252}.$ This implies $J=(4,7,9)$ and $\Phi(J)=\frac{1}{252}$. Then we consider the vertex $L$. For $L=(3,3,4,9),$ we get {the desired result analogous to the case for the vertex $I$.} For $L=(4,7,9)$, $M=(4,7,9)$ and $\Phi(M)=\frac{1}{252}.$ And in this case $N=(4,4,9)$ and $\Phi(N)=\frac{1}{9}>\frac{1}{12}.$ Hence $\Phi(A)+\Phi(B)+\Phi(D)+\Phi(I)+\Phi(J)+\Phi(N)>\frac{1}{12}$.

\item 
If $B=(3,3,4,k)$ for $7\leq k\leq 9$, see Figure \ref{case13}(b), then $\Phi(B)=\frac{1}{k}-\frac{1}{12}$. Since $\Phi(A)+\Phi(B)>\frac{1}{12}$ for $k=7,$ we only need to consider for $k=8,9$. Since the pattern of $C$ has to be $(4,7,k)$, $\Phi(C)=\frac{1}{k}-\frac{3}{28}$. The possible patterns of $D$ are $(3,3,3,k)$, $(3,3,4,k),$ $(3,8,24),$ and $(3,9,18)$.
If $D=(3,3,3,k)$, then $\Phi(A)+\Phi(B)+\Phi(C)+\Phi(D)>\frac{1}{12}$. For either $D=(3,8,24), k=8$ or $D=(3,9,18),k=9$, the pattern of $E$ is $(3,3,3,24)$ or $(3,3,3,18)$, which yields $\Phi(A)+\Phi(B)+\Phi(C)+\Phi(E)>\frac{1}{12}$.
For $D=(3,3,4,k)$, $\Phi(D)=\frac{1}{k}-\frac{1}{12}$. Hence if $k=8$, then $\Phi(A)+\Phi(B)+\Phi(C)+\Phi(D)>\frac{1}{12}$. So that it suffices to consider $k=9$. In this case the possible patterns of $F$ are $(4,7,9)$ and $(3,3,4,9)$. For $F=(3,3,4,9)$, $\Phi(F)=\frac{1}{36}$, which yields $\Phi(A)+\Phi(B)+\Phi(C)+\Phi(D)+\Phi(F)>\frac{1}{12}$. For $F=(4,7,9)$, $\Phi(F)=\frac{1}{252}.$ In this case $G=(4,7,9)$ and $\Phi(G)=\frac{1}{252}$. We can consider the vertex $H,$ similarly to $F.$ The only nontrivial pattern of $H$ is $(4,7,9)$, which implies $\Phi(H)=\frac{1}{252}$ and $\Phi(I)=\frac{1}{252}$. This yields $\Phi(J)\geq\frac{1}{36}$. Hence $\Phi(A)+\Phi(B)+\Phi(C)+\Phi(D)+\Phi(F)+\Phi(G)+\Phi(H)+\Phi(I)+\Phi(J)>\frac{1}{12}$.
 
\item 
 If $B=(4,8,8), $ see Figure \ref{case13}(c), then $\Phi(A)=\frac{1}{56}$ and $\Phi(B)=0$. And the pattern of $C$ has to be $(4,7,8)$, which implies $\Phi(C)=\frac{1}{56}$. The possible patterns of $D$ are $(4,7,l)$ and $(3,3,4,7)$.
For $D=(3,3,4,7)$, $\Phi(D)=\frac{5}{84}$, which yields $\Phi(A)+\Phi(C)+\Phi(D)>\frac{1}{12}$.
For $D=(4,7,l)$, we have $l\leq 8$ (otherwise $E$ has negative curvature). If $l=7$, then  $\Phi(A)+\Phi(C)+\Phi(D)+\Phi(E)>\frac{1}{12}$.
If $l=8$, then $\Phi(D)=\frac{1}{56}$,  $F=(4,7,8)$ and $\Phi(F)=\frac{1}{56}$. Note that $\Phi(H), \Phi(I)$ and $\Phi(G)$ are all at least $\frac{1}{252}$, and at least one of them has the curvature no less than $\frac{1}{56}$. So that  $\Phi(A)+\Phi(C)+\Phi(D)+\Phi(F)+\Phi(G)+\Phi(H)+\Phi(I)>\frac{1}{12}$.

\end{itemize}
\item[Case 14]
$A=(5,5,k)$ for $5\leq k \leq 9.$ In this case, $\Phi(A)=\frac{1}{k}-\frac{1}{10}$. Hence for $k=5$, $\Phi(A)=\frac{1}{12}+\frac{1}{60}>\frac{1}{12}$. So we only need to consider $k\geq 6$. Denote by $B$ and $C$ the neighbors of $A$ incident to the $k$-gon. The possible patterns of  $B$ are $(5,5,k)$, $(5,6,k)$ and $(3,3,5,k)$. 
  
 \begin{figure}[tb]
\begin{minipage}[b]{0.4\textwidth}
\centering
\includegraphics[width=4.2cm,height=4.2cm]{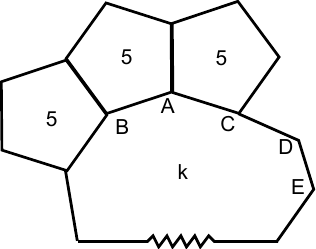}
(a)
\end{minipage}
\begin{minipage}[b]{0.4\textwidth}
\centering
\includegraphics[width=4.2cm,height=4.2cm]{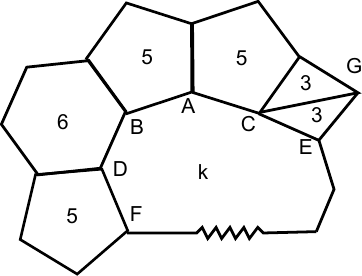}
(b)
\end{minipage}
\begin{minipage}[b]{0.4\textwidth}
\centering
\includegraphics[width=4.2cm,height=4.2cm]{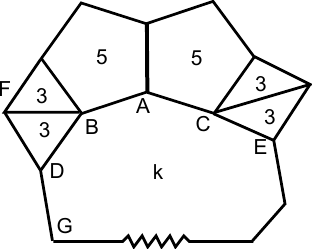}
(c)
\end{minipage}
\caption{\small Case 14.}
 \label{case14}
\end{figure}   
  
\begin{itemize}
\item 
If $B=(5,5,k)$ for {$6\leq k\leq 9,$} see Figure \ref{case14}(a), then $\Phi(B)=\frac{1}{k}-\frac{1}{10}$.
For $k\leq 7$, then $\Phi(A)+\Phi(B)>\frac{1}{12}$.
For $k=8,9$, then the only nontrivial pattern of $C$ is $(5,5,k)$, which implies $\Phi(C)=\frac{1}{k}-\frac{1}{10}$. Similarly, the curvature of other vertices on the $k$-polygon is equal to $\frac{1}{k}-\frac{1}{10}$, so that $\Phi(G)\geq k(\frac{1}{k}-\frac{1}{10})>\frac{1}{12}$.

\item 
If $B=(5,6,k)$ for $k=6,7,$ see Figure \ref{case14}(b), then $\Phi(B)=\frac{1}{k}-\frac{2}{15}$. Since $\Phi(A)+\Phi(B)>\frac{1}{12}$ for $k=6,$ it suffices to consider the case of $k=7$. In this case, $\Phi(A)=\frac{3}{70}$ and $\Phi(B)=\frac{1}{105}$. The only nontrivial pattern of $D$ is $(5,6,7)$, which implies $\Phi(D)=\frac{1}{105}$. The nontrivial patterns of $C$ are $(5,5,7),$ $(5,6,7)$ and $(3,3,5,7)$.
For $C=(5,5,7)$, $\Phi(C)=\frac{3}{70}$, which yields $\Phi(A)+\Phi(B)+\Phi(C)>\frac{1}{12}$. 
For $C=(5,6,7)$, $\Phi(C)=\frac{1}{105}$, and the pattern of $E$ has to be $(5,6,7)$, which implies $\Phi(E)=\frac{1}{105}$.
Noting that $\Phi(F)\geq\frac{1}{105},$ we have $\Phi(A)+\Phi(B)+\Phi(C)+\Phi(D)+\Phi(E)+\Phi(F)>\frac{1}{12}$. For $C=(3,3,5,7)$, $\Phi(C)=\frac{1}{105}$, and the pattern of $E$ is one of $(3,3,3,7)$, $(3,3,4,7),$ $(3,3,5,7)$ and $(3,7,42)$. If $E=(3,3,3,7)$, $(3,3,4,7),$ or $(3,3,5,7)$, then $\Phi(E)\geq\frac{1}{105}$. Hence $\Phi(A)+\Phi(B)+\Phi(C)+\Phi(D)+\Phi(E)+\Phi(F)>\frac{1}{12}$. If $E=(3,7,42)$, then the nontrivial pattern of $G$ is $(3,3,3,42)$, which implies $\Phi(G)=\frac{1}{42}$. Hence $\Phi(A)+\Phi(B)+\Phi(C)+\Phi(D)+\Phi(G)>\frac{1}{12}$.
 
 \item
If $B=(3,3,5,k),$ for {$6\leq k\leq7$} see Figure \ref{case14}(c), then $\Phi(B)=\frac{1}{k}-\frac{2}{15}$. Since $\Phi(A)+\Phi(B)>\frac{1}{12}$ for $k=6,$ we only need to consider the case of $k=7$. In this case, $\Phi(A)=\frac{3}{70}$ and $\Phi(B)=\frac{1}{105}$. The possible patterns of $C$ are $(5,5,7)$, $(5,6,7)$ and $(3,3,5,7)$. For $C=(5,5,7)$, $\Phi(C)=\frac{3}{70}$, which yields $\Phi(A)+\Phi(C)>\frac{1}{12}$. For $C=(5,6,7)$, {this is just the above case with the vertices $B,C$ exchanged.} For $C=(3,3,5,7)$, $\Phi(C)=\frac{1}{105}$. The possible patterns of $D$ and $E$ are $(3,3,3,7)$, $(3,3,4,7),$ $(3,3,5,7)$ and $(3,7,42)$. If $D=(3,3,3,7)$ or $(3,3,4,7)$, then $\Phi(D)\geq\frac{5}{84}$, which yields $\Phi(A)+\Phi(B)+\Phi(C)+\Phi(D)>\frac{1}{12}$. If $D=(3,7,42)$, then $F=(3,3,3,42)$ and $\Phi(F)=\frac{1}{42}$, which yields $\Phi(A)+\Phi(B)+\Phi(C)+\Phi(F)>\frac{1}{12}$. So the only nontrivial pattern of $D$ is $(3,3,5,7)$ and in this case $\Phi(D)=\frac{1}{105}.$ Similarly for the vertices $E$ and $G$, $\Phi(E),\Phi(G)=\frac{1}{105}$. Hence $\Phi(A)+\Phi(B)+\Phi(C)+\Phi(D)+\Phi(E)+\Phi(G)>\frac{1}{12}$.
  
 \end{itemize}
\item[Case 15]
$A=(5,6,k)$ for $6\leq k \leq 7.$ In this case, $\Phi(A)=\frac{1}{k}-\frac{2}{15}$. Denote by $B$ the neighbor of $A$ which is incident to the pentagon and the $k$-gon. The possible patterns of $B$ are $(5,6,k)$ and $(3,3,5,k)$.

 \begin{figure}[tb]
\begin{minipage}[b]{0.4\textwidth}
\centering
\includegraphics[width=4.2cm,height=4.2cm]{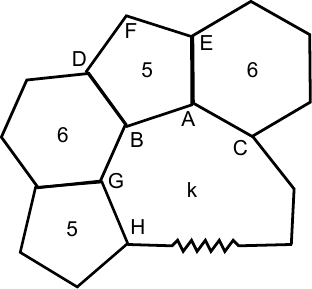}
(a)
\end{minipage}
\begin{minipage}[b]{0.4\textwidth}
\centering
\includegraphics[width=4.2cm,height=4.2cm]{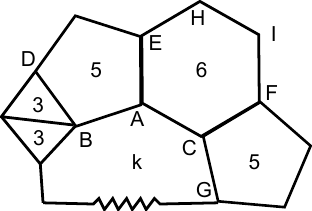}
(b)
\end{minipage}
\caption{\small Case 15.}
 \label{case15}
\end{figure}   
   
\begin{itemize}
\item 
If $B=(5,6,k)$, for {$6\leq k\leq 7$}, see Figure \ref{case15}(a), then $\Phi(B)=\frac{1}{k}-\frac{2}{15}$.
We consider the vertex $D$. The possible patterns of $D$ (and $E$) are 
$(5,6,6)$, $(5,6,7)$ and $(3,3,5,6)$. 
If $k=6$, then $\Phi(A)=\Phi(B)=\frac{1}{30}$ and $\Phi(D),\Phi(E)\geq\frac{1}{105}$.  
Hence $\Phi(A)+\Phi(B)+\Phi(D)+\Phi(E)>\frac{1}{12}$.
If $k= 7$, then $\Phi(A)=\Phi(B)=\frac{1}{105}$. And the pattern of $G$ or $C$ has to be $(5,6,7)$, which implies $\Phi(C)=\Phi(G)=\frac{1}{105}$ and $\Phi(H)\geq\frac{1}{105}$. Note that at least one of two vertices $D$ and $E$ is not of the pattern $(5,6,7)$ (otherwise, the vertex $F$ has negative curvature). Hence $\Phi(D)+\Phi(E)\geq\frac{1}{30}+\frac{1}{105}=\frac{3}{70}$.
So that $\Phi(A)+\Phi(B)+\Phi(C)+\Phi(D)+\Phi(E)+\Phi(G)+\Phi(H)>\frac{1}{12}$.

\item
If $B=(3,3,5,k)$, for {$6\leq k\leq 7$}, see Figure \ref{case15}(b), then for $k=6$, $\Phi(A)=\frac{1}{30}$,  $\Phi(B)=\frac{1}{30}$ and $\Phi(D), \Phi(E)\geq\frac{1}{105}$. Hence $\Phi(A)+\Phi(B)+\Phi(D)+\Phi(E)>\frac{1}{12}$.
For $k=7$, $\Phi(A)=\Phi(B)=\frac{1}{105}$. The pattern of $C$ has to be $(5,6,7)$, which implies $\Phi(C)=\frac{1}{105}$. Note that $\Phi(F),\Phi(D),\Phi(G)\geq\frac{1}{105}$. The possible patterns of $E$ are $(5,6,6)$, $(5,6,7)$ and $(3,3,5,6)$. If $E=(5,6,6)$ or $(3,3,5,6)$, then $\Phi(E)=\frac{1}{30}$, which yields $\Phi(A)+\Phi(B)+\Phi(C)+\Phi(D)+\Phi(E)+\Phi(F)+\Phi(G)>\frac{1}{12}$. If $E=(5,6,7)$, then $H=(5,6,7)$ and $\Phi(H)=\frac{1}{105}$. So that $\Phi(I)\geq\frac{1}{105}$. Hence $\Phi(A)+\Phi(B)+\Phi(C)+\Phi(D)+\Phi(E)+\Phi(F)+\Phi(G)+\Phi(H)+\Phi(I)>\frac{1}{12}$.
 
\end{itemize}

\item[Case 16]
$A=(3,3,3,k).$ In this case, $\Phi(A)=\frac{1}{k}$. For $k<12$, $\Phi(A)>\frac{1}{12}$. For $k=12$, we can construct a graph with total curvature $=\frac{1}{12}$, {see Figure \ref{fig1}}. So we only need to consider $k>12$. Denote by $B$ and $C$ the neighbors of $A$ incident to the $k$-gon.
The possible patterns of $B$ are $(3,3,3,k)$, $(3,7,42)$, $(3,8,24)$, $(3,9,18)$ and $(3,10,15)$.

 \begin{figure}[tb]
\begin{minipage}[b]{0.4\textwidth}
\centering
\includegraphics[width=4.2cm,height=4.2cm]{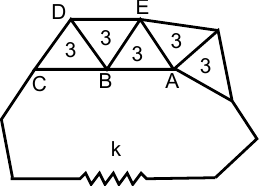}
(a)
\end{minipage}
\begin{minipage}[b]{0.4\textwidth}
\centering
\includegraphics[width=4.2cm,height=4.2cm]{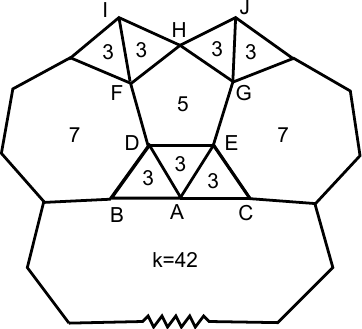}
(b)
\end{minipage}
\begin{minipage}[b]{0.4\textwidth}
\centering
\includegraphics[width=4.2cm,height=4.2cm]{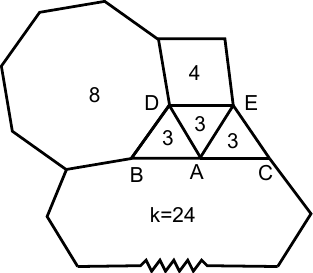}
(c)
\end{minipage}
\begin{minipage}[b]{0.4\textwidth}
\centering
\includegraphics[width=4.2cm,height=4.2cm]{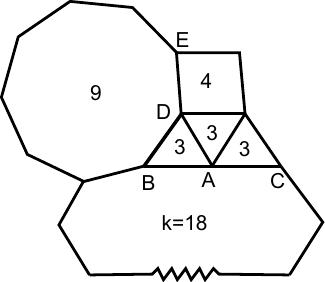}
(d)
\end{minipage}
\begin{minipage}[b]{0.4\textwidth}
\centering
\includegraphics[width=4.2cm,height=4.2cm]{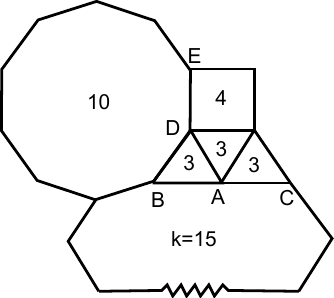}
(e)
\end{minipage}
\caption{\small Case 16.}
 \label{case16}
\end{figure}   

\begin{itemize}
\item 
If $B=(3,3,3,k)$, see Figure \ref{case16}(a), then $\Phi(A)=\Phi(B)=\frac{1}{k}$. For $k<24$, then $\Phi(A)+\Phi(B)>\frac{1}{12}$. We consider the case of $k=24$. If $\Phi(C)=0$, then $C=(3,8,24)$ and $\Phi(D)\geq\frac{1}{24}$. If $\Phi(C)\neq0$, then $\Phi(C)\geq\frac{1}{56}$. Hence we always have $\Phi(C)+\Phi(D)\geq\frac{1}{56}$.  So for $k=24$ we have $\Phi(A)+\Phi(B)+\Phi(C)+\Phi(D)>\frac{1}{12}$.
For $24<k<42$, the nontrivial pattern of other vertices on the $k$-gon is $(3,3,3,k).$ That is, the curvature of each vertex of the $k$-gon is equal to $\frac{1}{k}$, so that $\Phi(G)\geq1>\frac{1}{12}$.  
For $k=42$, the possible patterns of $C$ are $(3,3,3,42)$ and $(3,7,42)$. If $C=(3,7,42)$, then the possible patterns of $D$ are $(3,3,3,7)$, $(3,3,4,7)$ and $(3,3,5,7)$.
If $D=(3,3,3,7)$ or $(3,3,4,7)$, then $\Phi(D)\geq\frac{5}{84}$. Hence $\Phi(A)+\Phi(B)+\Phi(D)>\frac{1}{12}$. If $D=(3,3,5,7)$, then $\Phi(D)=\frac{1}{105}$ and the possible patterns of $E$ are $(3,3,3,5)$ and $(3,3,3,3,5)$, which implies $\Phi(E)\geq\frac{1}{30}.$ Hence $\Phi(A)+\Phi(B)+\Phi(D)+\Phi(E)>\frac{1}{12}$.
\item 
If $B=(3,7,42), $ see Figure \ref{case16}(b), then $\Phi(A)=\frac{1}{42}$ and $\Phi(B)=0$. We consider the vertex $C$. The possible patterns of $C$ are $(3,7,42)$ and $(3,3,3,42).$ 
For $C=(3,3,3,42),$ this is the same case as above (with $B,C$ exchanged). So we only need to consider $C=(3,7,42)$. In this case, the possible patterns of $D$ and $E$ are $(3,3,3,7)$, $(3,3,4,7)$ and $(3,3,5,7)$. However, if neither $D$ or $E$ is $(3,3,5,7)$, then $\Phi(A)+\Phi(D)+\Phi(E)>\frac{1}{12}$.
So the only nontrivial pattern of $D$ (or $E$) is $(3,3,5,7),$ which implies $\Phi(D)=\Phi(E)=\frac{1}{105}$. Note that $F=(3,3,5,7)$ and $G=(3,3,5,7)$, which implies $\Phi(F)=\Phi(G)=\frac{1}{105}$. Moreover, the only nontrivial pattern of $H$ is $(3,3,5,7)$, which yields $\Phi(H)\geq\frac{1}{105}$ and $\Phi(I),\Phi(J)\geq\frac{1}{105}$. Hence $\Phi(A)+\Phi(D)+\Phi(E)+\Phi(F)+\Phi(G)+\Phi(H)+\Phi(I)+\Phi(J)>\frac{1}{12}$.
\item 
If $B=(3,8,24),$ see Figure \ref{case16}(c), then 
$\Phi(A)=\frac{1}{24}$ and $\Phi(B)=0$. The pattern of $D$ has to be $(3,3,4,8)$ ({since if $D=(3,3,3,8)$, then $\Phi(D)=\frac{1}{8}>\frac{1}{12}$.}), which implies $\Phi(D)=\frac{1}{24}$. Moreover, it is easy to see that $\Phi(C)+\Phi(E)\geq\frac{1}{24}.$ 
Hence $\Phi(A)+\Phi(D)+\Phi(C)+\Phi(E)>\frac{1}{12}$.

\item 
If $B=(3,9,18), $ see Figure \ref{case16}(d), then 
$\Phi(A)=\frac{1}{18}$ and $\Phi(B)=0$. The only nontrivial case is $D=(3,3,4,9),$ which implies $\Phi(D)=\frac{1}{36}$. Since $\Phi(E)\geq\frac{1}{252}$, 
 $\Phi(A)+\Phi(D)+\Phi(E)>\frac{1}{12}$.

\item 
If $B=(3,10,15), $  see Figure \ref{case16}(e) then 
$\Phi(A)=\frac{1}{15}, \Phi(B)=0$. The only nontrivial case is $D=(3,3,4,10),$ which implies $\Phi(D)=\frac{1}{60}$.  Since $\Phi(E)\geq\frac{1}{60}$, 
$\Phi(A)+\Phi(D)+\Phi(E)>\frac{1}{12}$.
\end{itemize}

\item[Case 17]
$A=(3,3,4,k)$ for $4\leq k \leq 11.$ In this case, $\Phi(A)=\frac{1}{k}-\frac{1}{12}$. So that $\Phi(A)>\frac{1}{12}$ for $k=4,5$ and $\Phi(A)=\frac{1}{12}$ for $k=6.$ Moreover, we can construct graphs with total curvature $\frac{1}{12}$ for $k=6,11$, {see Figures \ref{fig3} and \ref{fig4}}. So we shall prove that for $7\leq k \leq 10$, $\Phi(G)>\frac{1}{12}$ (or it is reduced to the previous cases), and for $k=11$, $\Phi(G)\geq\frac{1}{12}$. Denote by $B$ a neighbor of $A$ which is incident to a triangle and the $k$-gon. The possible patterns of $B$ are $(3,3,4,k)$, $(3,3,5,7)$, $(3,7,42)$, $(3,8,24 )$, $(3,9,18)$ and $(3,10,15)$.

 \begin{figure}[tb]
\begin{minipage}[b]{0.4\textwidth}
\centering
\includegraphics[width=4.2cm,height=4.2cm]{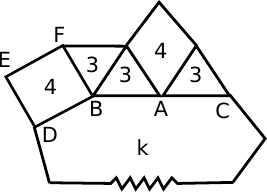}
(a)
\end{minipage}
\begin{minipage}[b]{0.4\textwidth}
\centering
\includegraphics[width=4.2cm,height=4.2cm]{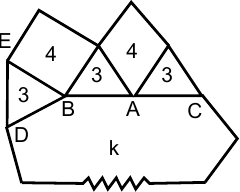}
(b)
\end{minipage}
\begin{minipage}[b]{0.4\textwidth}
\centering
\includegraphics[width=4.2cm,height=4.2cm]{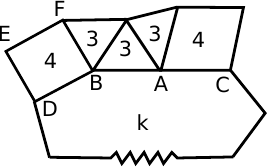}
(c)
\end{minipage}
\begin{minipage}[b]{0.4\textwidth}
\centering
\includegraphics[width=4.2cm,height=4.2cm]{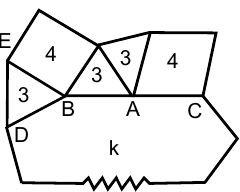}
(d)
\end{minipage}
\caption{\small Case 17: (a)-(d).}
 \label{case17 (a)-(d)}
\end{figure}

 \begin{figure}[tb]
\begin{minipage}[b]{0.4\textwidth}
\centering
\includegraphics[width=4.2cm,height=4.2cm]{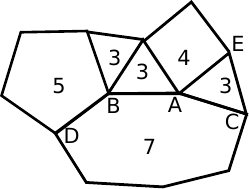}
(e)
\end{minipage}
\begin{minipage}[b]{0.4\textwidth}
\centering
\includegraphics[width=4.2cm,height=4.2cm]{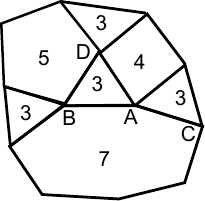}
(f)
\end{minipage}
\begin{minipage}[b]{0.4\textwidth}
\centering
\includegraphics[width=4.2cm,height=4.2cm]{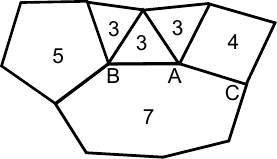}
(g)
\end{minipage}
\begin{minipage}[b]{0.4\textwidth}
\centering
\includegraphics[width=4.2cm,height=4.2cm]{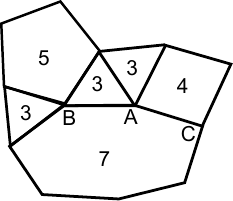}
(h)
\end{minipage}
\begin{minipage}[b]{0.4\textwidth}
\centering
\includegraphics[width=4.2cm,height=4.2cm]{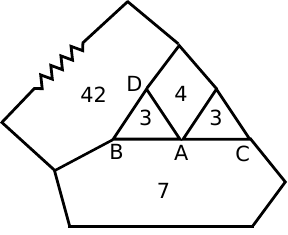}
(i)
\end{minipage}
\begin{minipage}[b]{0.4\textwidth}
\centering
\includegraphics[width=4.2cm,height=4.2cm]{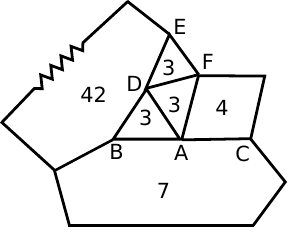}
(j)
\end{minipage}
\caption{\small Case 17: (e)-(j).}
 \label{case17 (e)-(j)}
\end{figure}

\begin{itemize}
\item  
If $B=(3,3,4,k)$, for {$7\leq k\leq11$} see Figure \ref{case17 (a)-(d)}(a-d), $\Phi(B)=\frac{1}{k}-\frac{1}{12}$. Since $\Phi(A)+\Phi(B)>\frac{1}{12}$ for $k=7$, we only need to consider $8\leq k\leq11$. So that the possible patterns of $D$ are $(3,3,4,k)$, $(3,8,24)$, $(3,9,18)$, $(3,10,15)$ and $(4,8,8)$. If $D=(3,8,24), (3,9,18), (3,10,15)$ or $(4,8,8),$ then $k= 8,9,10$ or $8,$ respectively. However, for $D=(3,8,24)$, $(3,9,18)$ or $(3,10,15)$, $\Phi(E)\neq0$ and hence $E$ is reduced to Cases 1-16 (with $A=E$ therein). For $D=(4,8,8)$, it is easy to see that $\Phi(E)+\Phi(F)\geq\frac{1}{56}$, which yields $\Phi(A)+\Phi(B)+\Phi(E)+\Phi(F)>\frac{1}{12}$. So the only nontrivial case is $D=(3,3,4,k)$, which implies $\Phi(D)=\frac{1}{k}-\frac{1}{12}$. The curvatures of other vertices of the $k$-polygon are equal to $\frac{1}{12},$ and hence $\Phi(G)\geq1-\frac{k}{12}$.

\item  
If $B=(3,3,5,7), $ then we have four subcases:
\begin{enumerate}
\item The first case is shown in Figure \ref{case17 (e)-(j)}(e), in which $\Phi(A)=\frac{5}{84}, $  $\Phi(B)=\frac{1}{105}$ and $\Phi(D)\geq\frac{1}{105}$. The nontrivial patterns of $C$ are $(3,3,4,7)$, $(3,3,5,7)$ and $(3,7,42)$.
For $C=(3,3,4,7)$, $\Phi(C)=\frac{5}{84}$, which yields $\Phi(A)+\Phi(C)>\frac{1}{12}$.
For $C=(3,3,5,7)$,  $\Phi(C)=\frac{1}{105}$, which implies $\Phi(A)+\Phi(B)+\Phi(C)+\Phi(D)>\frac{1}{12}$. For $C=(3,7,42),$ $E=(3,4,42)$, and hence $\Phi(A)+\Phi(E)>\frac{1}{12}$.

\item 
The second case is shown in Figure \ref{case17 (e)-(j)}(f), where $\Phi(A)=\frac{5}{84}$ and $\Phi(B)=\frac{1}{105}$. Since $D=(3,3,4,5)$, $\Phi(A)+\Phi(D)>\frac{1}{12}$.

\item 
The third case is shown in Figure \ref{case17 (e)-(j)}(g), where $\Phi(A)=\frac{5}{84}$ and $\Phi(B)=\frac{1}{105}.$ The only nontrivial case is $C=(3,3,4,7)$, which implies $\Phi(C)=\frac{5}{84}$. Hence $\Phi(A)+\Phi(C)>\frac{1}{12}$.

\item The fourth case is shown in Figure \ref{case17 (e)-(j)} (h). The proof is similar and hence omitted.
\end{enumerate}

\item
If $B=(3,7,42), $ then we have two subcases:
\begin{enumerate}
 \item  The first case is shown in Figure \ref{case17 (e)-(j)}(i). Since $D=(3,4,42)$, $\Phi(A)+\Phi(D)>\frac{1}{12}$.
 
 \item  The second case is shown in Figure \ref{case17 (e)-(j)}(j). In this case, $\Phi(A)=\frac{5}{84}$. Since $D=(3,3,3,42)$, $\Phi(D)=\frac{1}{42}$. Clearly, $\Phi(E)+\Phi(F)\geq\frac{1}{42}$ and $\Phi(A)+\Phi(D)+\Phi(E)+\Phi(F)>\frac{1}{12}$. 
 
\end{enumerate}
The same arguments work for $B=(3,8,24)$, $(3,9,18)$ and $(3,10,15)$ and hence are omitted here.
 \end{itemize}

\item[Case 18]
$A=(3,3,5,k)$ for $5\leq k \leq 7.$ In this case, $\Phi(A)=\frac{1}{k}-\frac{2}{15}$. Denote by $B$ a neighbor of $A$ which is incident to a triangle and the $k$-gon. The possible patterns of $B $ are $(3,3,5,k)$, $(3,4,4,5)$, $(3,3,3,3,5)$, $(3,7,42 )$, $(3,3,6,6)$, $(3,4,4,6)$ and $(3,3,3,3,6)$.

 \begin{figure}[tb]
\begin{minipage}[b]{0.4\textwidth}
\centering
\includegraphics[width=4.2cm,height=4.2cm]{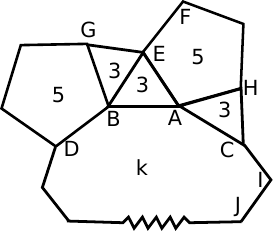}
(a)
\end{minipage}
\begin{minipage}[b]{0.4\textwidth}
\centering
\includegraphics[width=4.2cm,height=4.2cm]{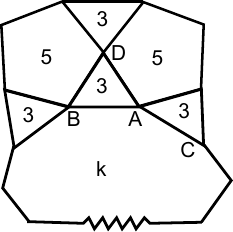}
(b)
\end{minipage}
\begin{minipage}[b]{0.4\textwidth}
\centering
\includegraphics[width=4.2cm,height=4.2cm]{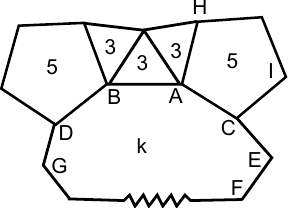}
(c)
\end{minipage}
\begin{minipage}[b]{0.4\textwidth}
\centering
\includegraphics[width=4.2cm,height=4.2cm]{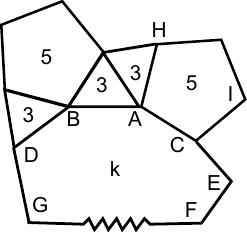}
(d)
\end{minipage}
\begin{minipage}[b]{0.4\textwidth}
\centering
\includegraphics[width=4.2cm,height=4.2cm]{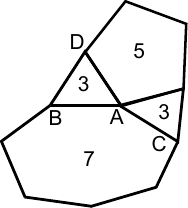}
(e)
\end{minipage}
\begin{minipage}[b]{0.4\textwidth}
\centering
\includegraphics[width=4.2cm,height=4.2cm]{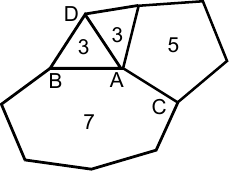}
(f)
\end{minipage}
\caption{\small Case 18: (a)-(f).}
 \label{case18 (a)-(f)}
\end{figure}

 \begin{figure}[tb]
\begin{minipage}[b]{0.4\textwidth}
\centering
\includegraphics[width=4.2cm,height=4.2cm]{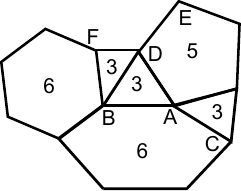}
(g)
\end{minipage}
\begin{minipage}[b]{0.4\textwidth}
\centering
\includegraphics[width=4.2cm,height=4.2cm]{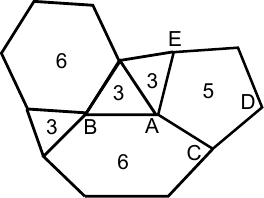}
(h)
\end{minipage}
\begin{minipage}[b]{0.4\textwidth}
\centering
\includegraphics[width=4.2cm,height=4.2cm]{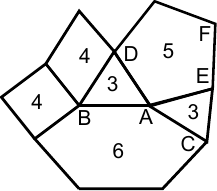}
(i)
\end{minipage}
\begin{minipage}[b]{0.4\textwidth}
\centering
\includegraphics[width=4.2cm,height=4.2cm]{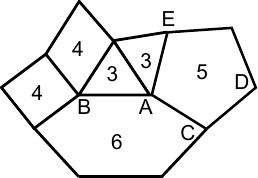}
(j)
\end{minipage}
\begin{minipage}[b]{0.4\textwidth}
\centering
\includegraphics[width=4.2cm,height=4.2cm]{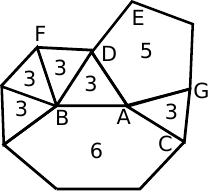}
(k)
\end{minipage}
\begin{minipage}[b]{0.4\textwidth}
\centering
\includegraphics[width=4.2cm,height=4.2cm]{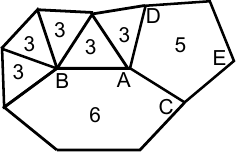}
(l)
\end{minipage}
\caption{\small Case 18: (g)-(l).}
 \label{case18 (g)-(l)}
\end{figure}

\begin{itemize}
\item  
If $B=(3,3,5,k), $ then $\Phi(B)=\frac{1}{k}-\frac{2}{15}$. Since $\Phi(A)+\Phi(B)>\frac{1}{12}$ for $k=5,$ it suffices to consider $k=6,7$.
\begin{enumerate}
\item 
The first case is shown in Figure \ref{case18 (a)-(f)} (a). The only nontrivial case is $D=(3,3,5,k)$, which implies
$\Phi(D)=\frac{1}{k}-\frac{2}{15}$. Hence $\Phi(A)+\Phi(B)+\Phi(D)>\frac{1}{12}$ for $k=6$. So we only need to consider $k=7$.
In this case, we consider the vertex $E$, the nontrivial patterns of $E$ are $(3,3,5,5)$, $(3,3,5,6)$ and $(3,3,5,7)$.
For $E=(3,3,5,5)$, $\Phi(E)=\frac{1}{15}$, which yields  $\Phi(A)+\Phi(B)+\Phi(D)+\Phi(E)>\frac{1}{12}$. For $E=(3,3,5,6)$, $\Phi(E)=\frac{1}{30}$. Then the pattern of $F$ is $(3,3,5,6)$, which implies $\Phi(F)=\frac{1}{30}$. Hence $\Phi(A)+\Phi(B)+\Phi(D)+\Phi(E)+\Phi(F)>\frac{1}{12}$. For $E=(3,3,5,7)$, $\Phi(E)=\Phi(F)=\Phi(G)=\frac{1}{105}$. Now we consider the vertex $C$. The possible patterns of $C$ are $(3,7,42)$ and $(3,3,5,7)$. If $C=(3,7,42)$, then $H=(3,5,42)$, which yields $\Phi(A)+\Phi(B)+\Phi(D)+\Phi(H)>\frac{1}{12}$. Hence the only nontrivial case is $C=(3,3,5,7)$, which implies $\Phi(C)=\frac{1}{105}$. Similarly, the nontrivial pattern of $I$ (or $J$) is also $(3,3,5,7)$, which implies $\Phi(I)=\Phi(J)=\frac{1}{105}$. Hence
$\Phi(A)+\Phi(B)+\Phi(C)+\Phi(D)+\Phi(E)+\Phi(F)+\Phi(G)+\Phi(I)+\Phi(J)>\frac{1}{12}$.

\item  
The second case is shown in Figure \ref{case18 (a)-(f)} (b). Since $D=(3,3,5,5)$, we have $\Phi(D)=\frac{1}{15}$, which yields $\Phi(A)+\Phi(B)+\Phi(D)>\frac{1}{12}$. 

\item
The third case is shown in Figure \ref{case18 (a)-(f)} (c). Since the nontrivial pattern of $D$ is $(3,3,5,k)$, $\Phi(D)=\frac{1}{k}-\frac{2}{15}$. For $k=6$, $\Phi(A)+\Phi(B)+\Phi(D)>\frac{1}{12}$. So we only need to consider $k=7$. In this case, $C=(3,3,5,7)$ and $\Phi(C)=\frac{1}{105}$. Similar to the subcase (1), we consider the vertices $E,F,G$ sequently. The nontrivial pattern of $E,F$ or $G$ is $(3,3,5,7)$, which implies $\Phi(E)=\Phi(F)=\Phi(G)=\frac{1}{105}$. Noting that $\Phi(H),\Phi(I)\geq\frac{1}{105}$, we have
$\Phi(A)+\Phi(B)+\Phi(C)+\Phi(D)+\Phi(E)+\Phi(F)+\Phi(G)+\Phi(H)+\Phi(I)>\frac{1}{12}$.

\item
The fourth case is shown in Figure \ref{case18 (a)-(f)} (d). The proof is similar to the above case, and hence is omitted.  
\end{enumerate}

\item  
If $B=(3,4,4,5)$ or $(3,3,3,3,5)$, then $\Phi(A)=\frac{1}{15}$ and $\Phi(B)=\frac{1}{30}$, which yields $\Phi(A)+\Phi(B)>\frac{1}{12}$.

\item
If $B=(3,7,42)$, see Figure \ref{case18 (a)-(f)} (e, f), then in both cases, the pattern of $D$ is $(3,3,42)$, $(3,5,42)$ or $(3,3,3,42)$, and it reduces to Case 1, Case 3 or Case 16, respectively (with $A=D$). So we get the desired conclusion.

\item If $B=(3,3,6,6)$, then there are four subcases. {However the proofs are similar, so we consider the following two typical cases:} Since their proofs are similar, we only consider the following typical cases:
\begin{enumerate}
\item
The first case is shown in Figure \ref{case18 (g)-(l)}(g) where $\Phi(A)=\frac{1}{30}$ and $\Phi(B)=0$. We consider the vertex $D$. For $D=(3,3,5,5),$ $\Phi(D)=\frac{1}{15}$, which yields $\Phi(A)+\Phi(D)>\frac{1}{12}$.
For $D=(3,3,5,6)$, $\Phi(D)=\frac{1}{30}$ and $\Phi(E)=\frac{1}{30}$, which yields $\Phi(A)+\Phi(D)+\Phi(E)>\frac{1}{12}$. For $D=(3,3,5,7)$, $F=(3,6,7)$, and hence $\Phi(A)+\Phi(D)+\Phi(F)>\frac{1}{12}$.
 \item
The second case is shown in Figure \ref{case18 (g)-(l)}(h) where $\Phi(A)=\frac{1}{30}$ and $\Phi(B)=0$. Note that the pattern of $C$ has to be $(3,3,5,6)$, which implies $\Phi(C)=\frac{1}{30}$. And it is obvious that $\Phi(D),\Phi(E)\geq\frac{1}{105}$, which yields $\Phi(A)+\Phi(C)+\Phi(D)+\Phi(E)>\frac{1}{12}$.

 \end{enumerate}
\item If $B=(3,4,4,6),$ {then we have two subcases:}

\begin{enumerate}
\item
The first one is shown in Figure \ref{case18 (g)-(l)} (i) where $\Phi(A)=\frac{1}{30}$ and  $\Phi(B)=0$. We consider the vertex $D$. The nontrivial patterns of $D$ are $(3,3,4,5)$ and $(3,4,4,5)$. For $D=(3,3,4,5)$, $\Phi(A)+\Phi(D)>\frac{1}{12}$. For $D=(3,4,4,5)$, $\Phi(D)=\frac{1}{30}$. Then we consider the vertex $E$. If $E=(3,3,5,5)$ or $(3,3,5,6)$, then $\Phi(E)\geq\frac{1}{30}$, which yields $\Phi(A)+\Phi(D)+\Phi(E)>\frac{1}{12}$. If $E=(3,3,5,7)$, then $\Phi(E)=\frac{1}{105}$ and $\Phi(F)\geq\frac{1}{105}$, which yields $\Phi(A)+\Phi(D)+\Phi(E)+\Phi(F)>\frac{1}{12}$.
 
\item
The other is shown in Figure \ref{case18 (g)-(l)} (j) where $\Phi(A)=\frac{1}{30}$ and $\Phi(B)=0$. Note that the pattern $C$ is $(3,3,5,6)$, which implies $\Phi(C)=\frac{1}{30}$. It is obvious that $\Phi(D),\Phi(E)\geq\frac{1}{105}$, and hence $\Phi(A)+\Phi(C)+\Phi(D)+\Phi(E)>\frac{1}{12}$.
\end{enumerate} 
  
\item 
If $B=(3,3,3,3,6),$ then there exist the following two subcases:   

\begin{enumerate}
\item The first one is shown in Figure \ref{case18 (g)-(l)} (k) where $\Phi(A)=\frac{1}{30}$ and $\Phi(B)=0$. We consider the vertex $D$. The possible patterns of $D$ are $(3,3,5,5)$, $(3,3,5,6)$, $(3,3,5,7)$ and $(3,3,3,3,5)$. For $D=(3,3,5,5)$, $\Phi(D)=\frac{1}{15}$, which yields $\Phi(A)+\Phi(D)>\frac{1}{12}$. For $D=(3,3,5,6)$, $\Phi(D)=\frac{1}{30}$ and $\Phi(E)\geq\frac{1}{30}$, which yields $\Phi(A)+\Phi(D)+\Phi(E)>\frac{1}{12}$. For $D=(3,3,3,3,5)$, $\Phi(D)=\frac{1}{30},$ $\Phi(E)\geq\frac{1}{105}$ and $\Phi(G)\geq\frac{1}{105}$. Hence $\Phi(A)+\Phi(D)+\Phi(E)+\Phi(G)>\frac{1}{12}$. For $D=(3,3,5,7)$, the nontrivial pattern of $F$ is $(3,3,5,7)$. In this case, we can regard the vertex $D$ as the vertex $A$ of Figure \ref{case18 (a)-(f)} (c,d,f) to get the desired result.
\item The other is shown in Figure \ref{case18 (g)-(l)} (l) where $\Phi(A)=\frac{1}{30}$ and $\Phi(B)=0$. Note that the nontrivial case is $C=(3,3,5,6)$, which implies $\Phi(C)=\frac{1}{30}$ and $\Phi(D),\Phi(E)\geq\frac{1}{105}$. Hence $\Phi(A)+\Phi(C)+\Phi(D)+\Phi(E)>\frac{1}{12}$.
\end{enumerate}
\end{itemize}

\item[Case 19]
$A=(3,4,4,k)$ for $4\leq k \leq 5.$ For $A= (3,4,4,4)$, {we can construct  a graph with total curvature $\frac{1}{12}$, see Figure \ref{fig2}}. So we only need to consider $k=5$ which is divided into two subcases:
 \begin{figure}[tb]
\begin{minipage}[b]{0.4\textwidth}
\centering
\includegraphics[width=4.2cm,height=4.2cm]{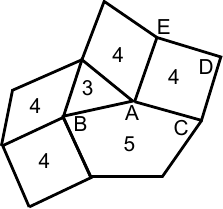}
(a)
\end{minipage}
\begin{minipage}[b]{0.4\textwidth}
\centering
\includegraphics[width=4.2cm,height=4.2cm]{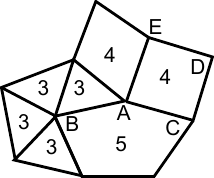}
(b)
\end{minipage}
\begin{minipage}[b]{0.4\textwidth}
\centering
\includegraphics[width=4.2cm,height=4.2cm]{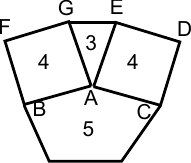}
(c)
\end{minipage}
\caption{\small Case 19.}
 \label{case19}
\end{figure}   

\begin{itemize}
\item 
The first case is shown in Figure \ref{case19}(a,b). The possible patterns of $B $ are $(3,4,4,5)$ and $(3,3,3,3,5)$.

If $B=(3,4,4,5),$ see Figure \ref{case19} (a), then $\Phi(A)=\frac{1}{30}$ and $\Phi(B)=\frac{1}{30}$. The possible patterns of $C$ are $(3,4,4,5)$ and $(4,5,20)$. For $C=(3,4,4,5)$, $\Phi(C)=\frac{1}{30}$, which yields $\Phi(A)+\Phi(B)+\Phi(C)>\frac{1}{12}$. For $C=(4,5,20)$, $D=(4,5,20)$ and $E=(3,4,4,5)$, which implies $\Phi(E)=\frac{1}{30}$. Hence $\Phi(A)+\Phi(B)+\Phi(E)>\frac{1}{12}$.

If $B=(3,3,3,3,5),$ see Figure \ref{case19} (b), then $\Phi(A)=\frac{1}{30}, $ and $\Phi(B)=\frac{1}{30}$. Similar to the above case, considering the vertices $C, D$ and $E$, we can get the desired conclusion.
 
 \item 
The second case is shown in Figure \ref{case19} (c).   The possible patterns of $B$ and $C$ are $(3,4,4,5)$ and $(4,5,20)$.

If the patterns of $B$ and $C$ are $(3,4,4,5)$,
then $\Phi(A)=\Phi(B)=\Phi(C)=\frac{1}{30}$, which yields $\Phi(A)+\Phi(B)+\Phi(C)>\frac{1}{12}$. If $B=(3,4,4,5)$ and $C=(4,5,20)$, then $\Phi(B)=\frac{1}{30}$ and the pattern of $D$ is $(4,5,20)$. Hence the pattern of $E$ has to be $(3,4,4,5)$ and $\Phi(E)=\frac{1}{30}$. Hence $\Phi(A)+\Phi(B)+\Phi(E)>\frac{1}{12}$. If the patterns of $B$ and $C$ are $(4,5,20)$, then $\Phi(A)=\frac{1}{30}$ and $\Phi(B)=\Phi(C)=0$. The patterns of $D$ and $F$ are $(4,5,20)$. Hence the nontrivial patterns of $E$ and $G$ are $(3,4,4,5),$ which implies $\Phi(E)=\Phi(G)=\frac{1}{30}$. Hence $\Phi(A)+\Phi(E)+\Phi(G)>\frac{1}{12}$.
\end{itemize}

\item[Case 20]
$A=(3,3,3,3,k)$ for $3\leq k\leq 5.$

 \begin{figure}[tb]
\begin{minipage}[b]{0.4\textwidth}
\centering
\includegraphics[width=4.2cm,height=4.2cm]{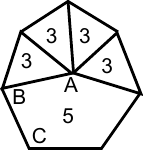}
\end{minipage}
\caption{\small Case 20.}
 \label{case20}
\end{figure}   

For $k=3$,  $\Phi(A)=\frac{1}{6}>\frac{1}{12}$. For $k=4$, {we can construct a graph with total curvature $\frac{1}{12}$, see Figure \ref{fig3}}. So we only need to consider $k=5.$ Denote by $B$ and $C$ neighbors of $A$ incident to the pentagon. The possible pattern of $B$ is $(3,3,3,3,5).$ In this case,
$\Phi(A)=\frac{1}{30}$ and $\Phi(B)=\frac{1}{30}$. Similarly, the possible pattern of $C$ is also $(3,3,3,3,5)$, which implies $\Phi(C)=\frac{1}{30}$. Hence $\Phi(A)+\Phi(B)+\Phi(C)>\frac{1}{12}$.
 \end{description}%
By combining all above cases, we prove the first part of the theorem.

In the above case-by-case proof, one can figure out that if the total curvature of $G$ attains the first gap of the total curvature, i.e. $\frac{1}{12}$, then the polygonal surface $S(G)$ has only two types of metric structures, $(a)$ or $(b)$, as indicated in the Theorem~\ref{mainthm}. This proves the second (rigidity) part of the theorem.
\end{proof}


\section{Examples}\label{s:examples}
\begin{figure}[htbp]
\begin{center}
\includegraphics[width=0.8\linewidth]{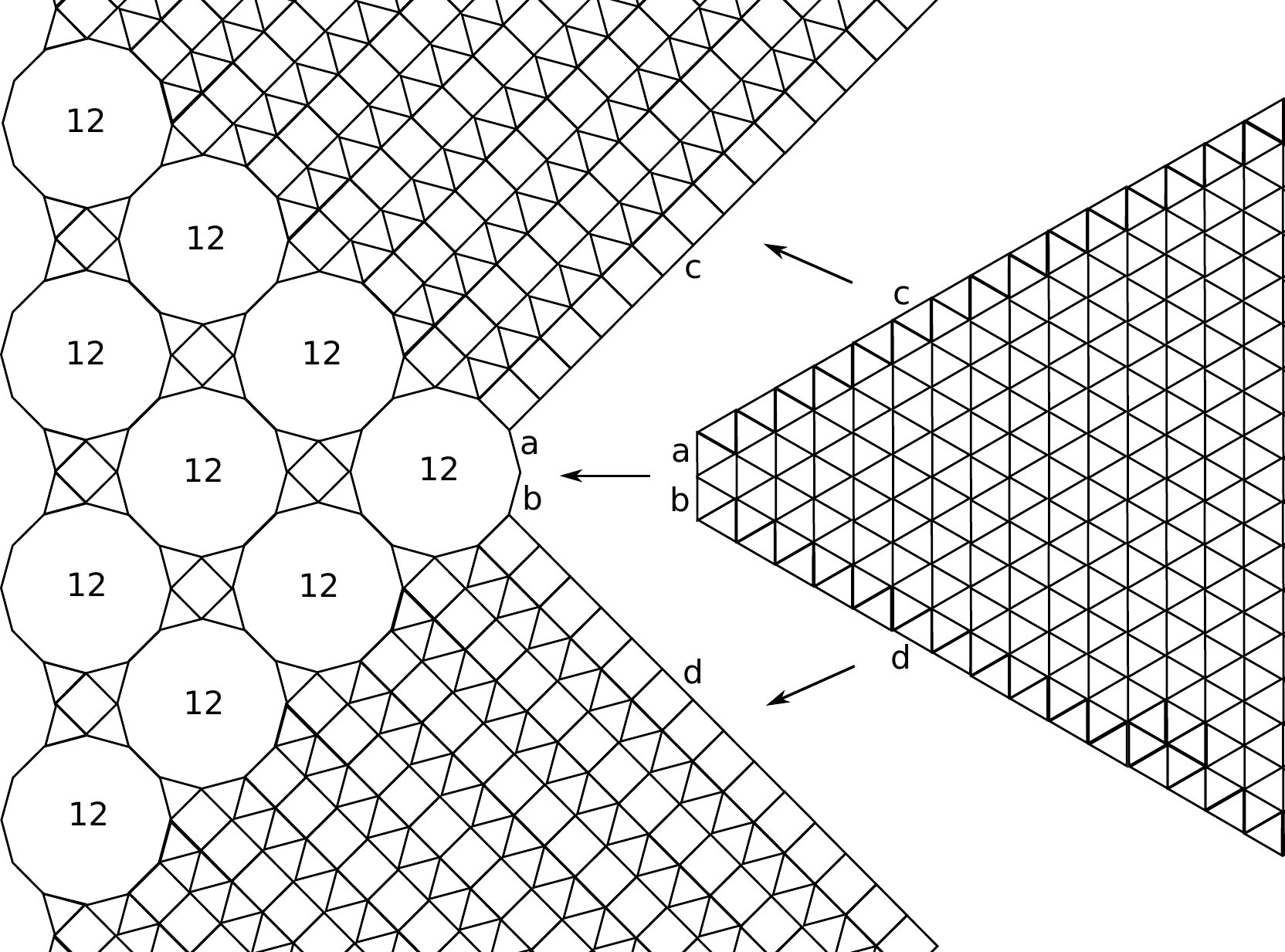}
\caption{\small A graph in $\mathcal{PC}_{\geq0}$ with total curvature $\frac{1}{12}$ which has a single vertex with non-vanishing curvature and of patterns $(3,6,12)$, $(3,3,3,12)$, $(3,3,4,6)$ and $(3,3,3,3,4)$.}
\label{fig1}
\end{center}
\end{figure}

\begin{figure}[htbp]
\begin{center}
\includegraphics[width=0.6\linewidth]{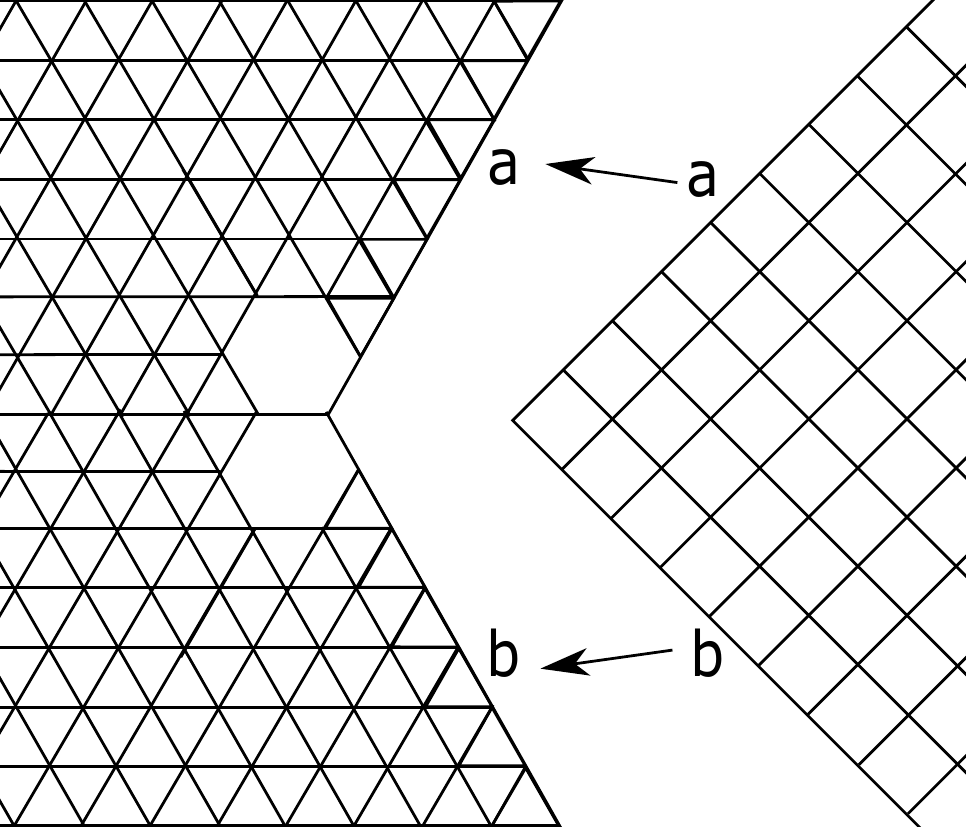}
\caption{\small A graph in $\mathcal{PC}_{\geq0}$ with total curvature $\frac{1}{12}$ which has a single vertex with non-vanishing curvature and of patterns $(4,6,6)$, $(3,3,4,6)$ and $(3,3,3,3,4)$.}
\label{fig3}
\end{center}
\end{figure}

\begin{figure}[htbp]
\begin{center}
\includegraphics[width=0.7\linewidth]{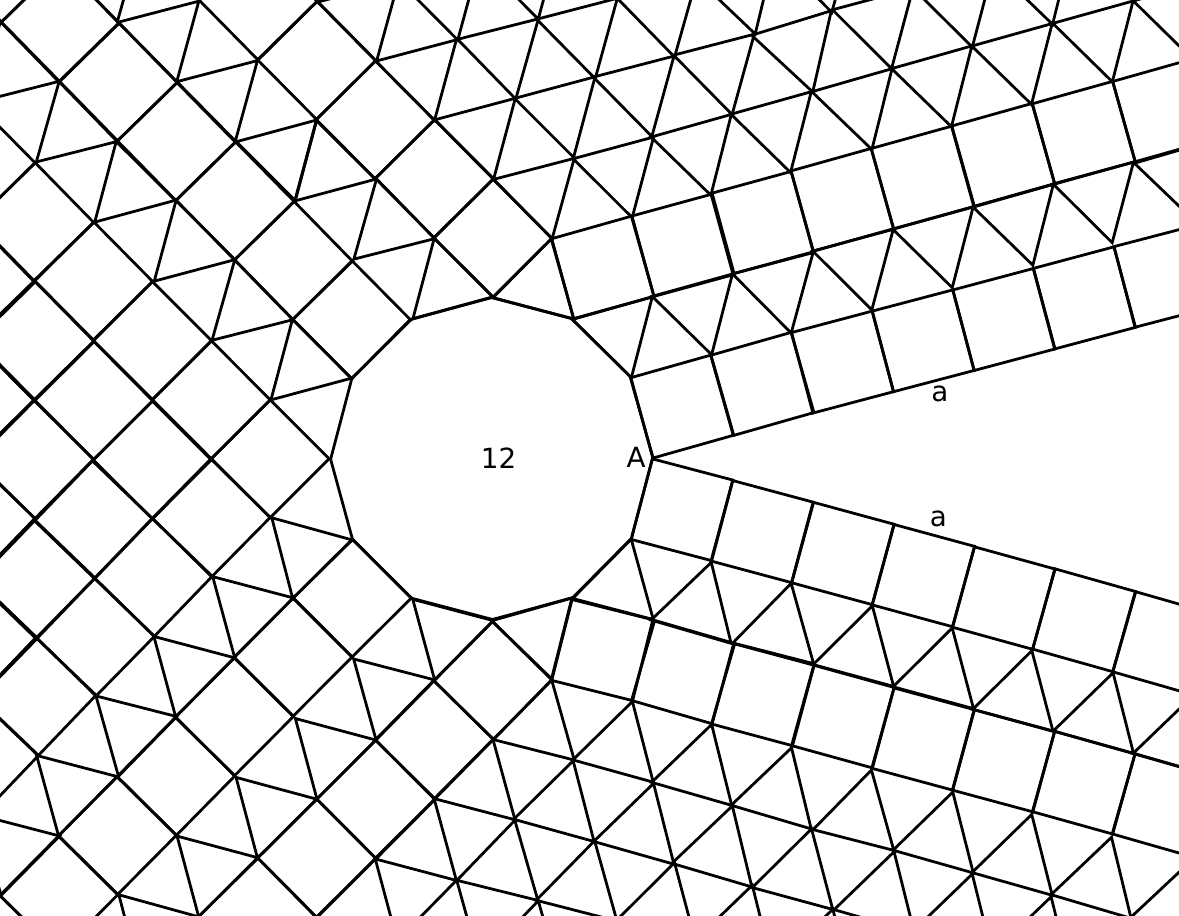}
\caption{\small A graph in $\mathcal{PC}_{\geq0}$ with total curvature $\frac{1}{12}$ which has a single vertex with non-vanishing curvature and of patterns $(4,4,12)$ and $(3,4,4,4).$}
\label{4-4-12}
\end{center}
\end{figure}

\begin{figure}[htbp]
\begin{center}
\includegraphics[width=0.75\linewidth]{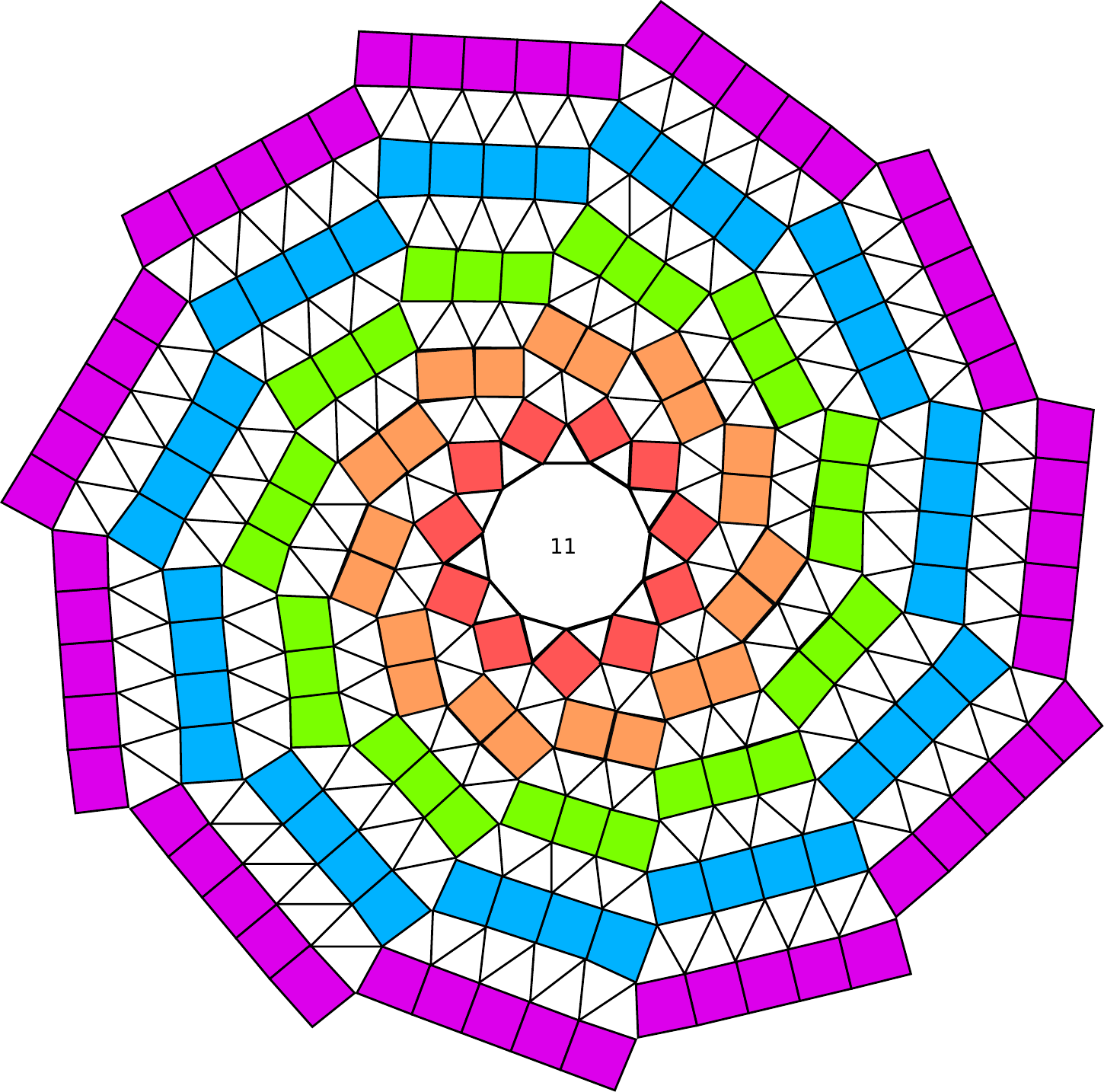}
\caption{\small A graph in $\mathcal{PC}_{\geq0}$ with total curvature $\frac{1}{12}$ which has eleven vertices (the vertices of the hendecagon) with non-vanishing curvature and of pattern $(3,3,4,11).$}
\label{fig4}
\end{center}
\end{figure}

\begin{figure}[htbp]
\begin{center}
\includegraphics[width=0.8\linewidth]{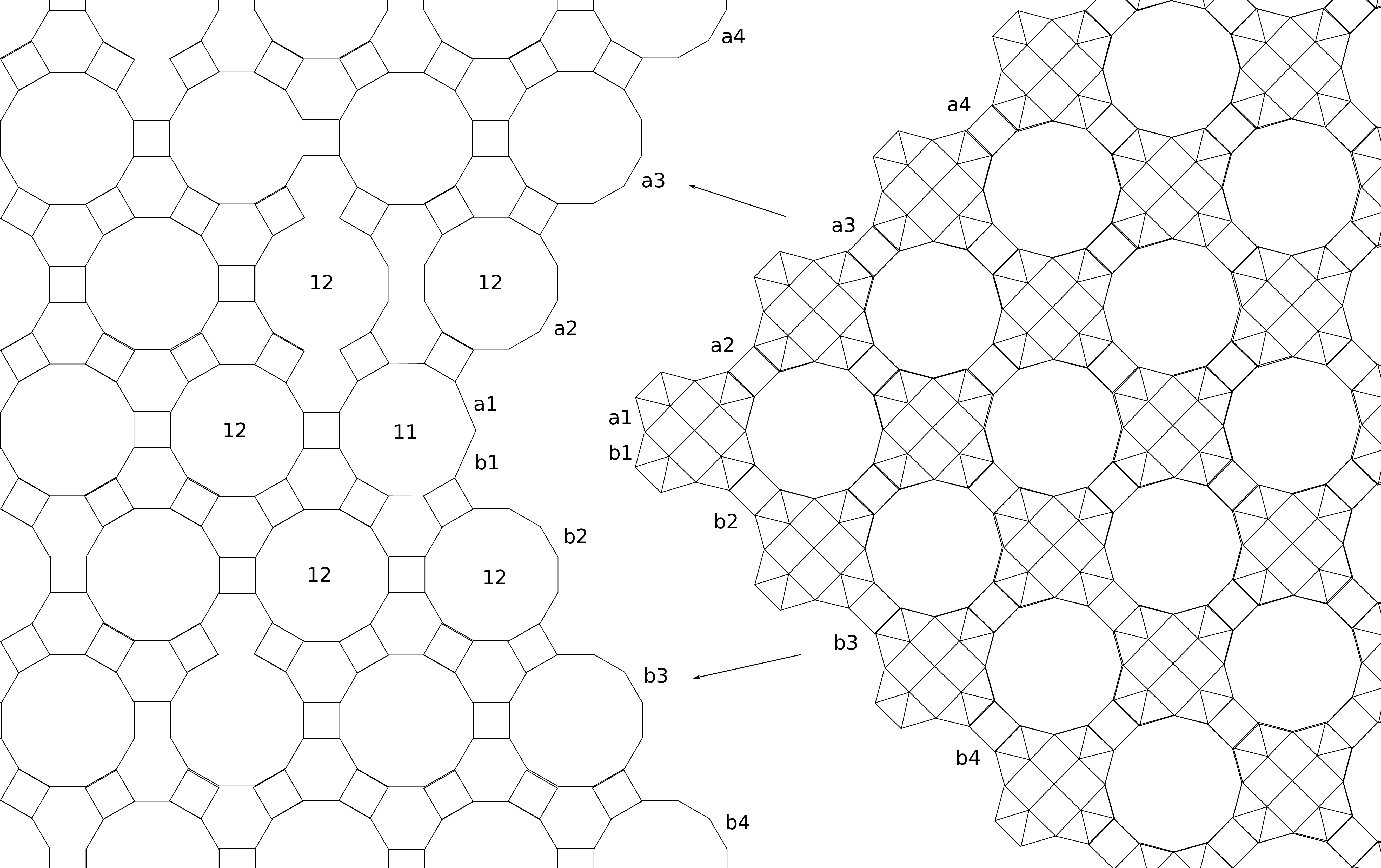}
\caption{\small A graph in $\mathcal{PC}_{\geq0}$ with total curvature $\frac{1}{12}$ which has eleven vertices (the vertices of the hendecagon) with non-vanishing curvature and of patterns $(4,6,11)$ and $(3,3,4,11).$}
\label{fig6}
\end{center}
\end{figure}

In this section, we construct some examples to show the sharpness of the main result, Theorem~\ref{mainthm}, 
see Figures~\ref{fig1}-\ref{fig6}.


\section{Applications}
\subsection{Connected components of induced graphs on $T(G)$}
Let $G=(V,E,F)$ be an infinite semiplanar graph with nonnegative curvature and $T(G)$ be the set of vertices with non-vanishing curvature. We are interested in the structure of $T(G),$ which, as a induced subgraph, is usually not connected. As a byproduct of the proof of Theorem~\ref{mainthm}, we give the upper bound of the number of connected components for the induced subgraph on $T(G).$ 
\begin{corollary}\label{coro:concomp} For an infinite semiplanar graph $G$ with nonnegative curvature, the induced subgraph on $T(G)$ has at most $14$ connected components. 
\end{corollary}
\begin{proof}

As in the proof of Theorem \ref{mainthm}, starting from 
a vertex $A\in T(G),$ we find several nearby vertices in $T(G)$ such that the sum of their curvatures is at least $\frac{1}{12}.$ For any $A\in T(G),$ we denote by $C_A$ the connected component of the induced subgraph on $T(G)$ which contains $A,$ and by $\Phi(C_A):=\sum_{x\in C_A} \Phi(x)$ the sum of the curvatures of vertices in $C_A.$ By checking the proof of Theorem \ref{mainthm}, we have the following:
\begin{enumerate}[(a)]
\item In a subcase of Case 8, see Figure \ref{case8: (f)-(i).}(f), $C_A$ consists of only two vertices of pattern $(3,10,10)$ and $\Phi(C_A)=\frac{1}{15}.$
\item In other cases, the vertices with non-vanishing curvature we found are connected and hence $\Phi(C_A)\geq \frac{1}{12}.$
\end{enumerate}
This yields that $\Phi(C_A)\geq \frac{1}{15}$ for any $A\in T(G).$ Since $\Phi(G)\leq 1,$ we obtain that the number of connected components of the induced subgraph on $T(G)$ is at most 15.

To prove the theorem, it suffices to exclude the case that the induced subgraph on $T(G)$ has 15 connected components. Suppose it is the case, then each connected component of the induced subgraph on $T(G)$ is in the case $(a)$ above. Let $A$ and $B$ be of pattern $(3,10,10)$ as in Figure \ref{15-component} such that $\Phi(C_A)=\frac{1}{15}.$ Then, one can show that the vertices $C,D,E,F,G$ and $H$ have vanishing curvature. Hence the vertex $I$ has non-vanishing curvature and is of pattern $(3,10,10).$ This yields that the pattern of $J$ is also $(3,10,10)$. Similar argument works for $K$ and $L$. Then by the combinatorial restriction, one can show that $M\in T(G)$. Hence $I,J,K,L$ and $M$ are in $T(G)$ and connected, which implies that $\Phi(C_I)>\frac{1}{15}$. This yields a contradiction and proves the corollary.

 \begin{figure}[tb]
\begin{minipage}[b]{0.3\textwidth}
\centering
\includegraphics[width=4.9cm,height=5.6cm]{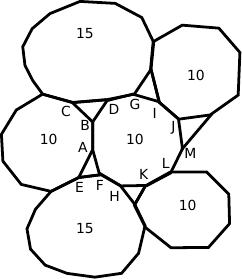}
\end{minipage}
\caption{\small Proof of Corollary~\ref{coro:concomp}}
 \label{15-component}
\end{figure}

\end{proof}

\begin{figure}[htbp]
\begin{center}
\includegraphics[width=0.6\linewidth]{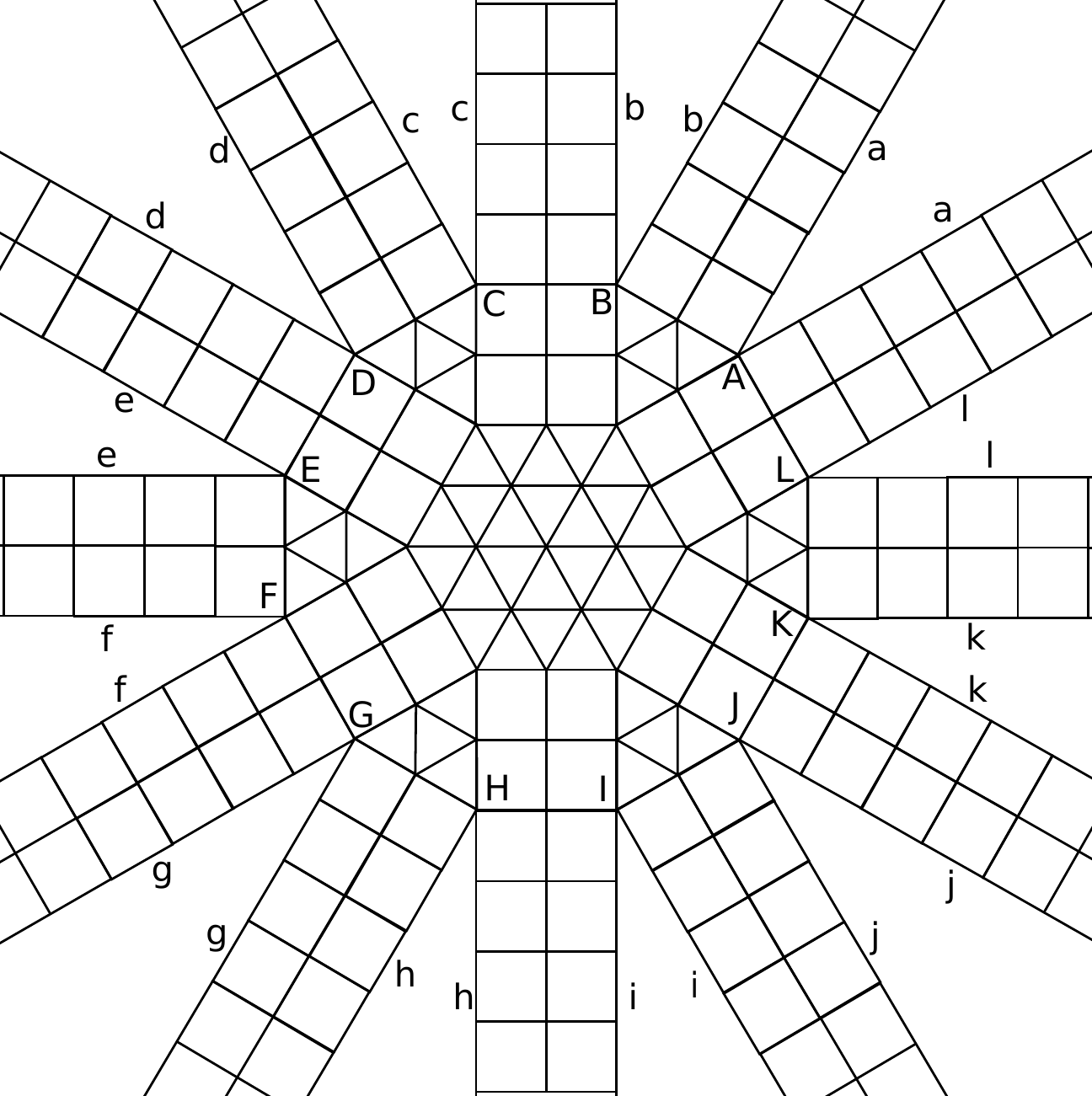}
\caption{\small The induced subgraph on $T(G)$ of a graph $G$ in $\NNG$ has 12 connected components.}
\label{12-component}
\end{center}
\end{figure}

Moreover, we construct a semiplanar graph $G$ whose induced subgraph on $T(G)$ has $12$ connected components, see Figure \ref{12-component}. Based on this, we propose the following conjecture.
\begin{conjecture} Let $G=(V,E,F)$ be a semiplanar graph with nonnegative curvature. The number of connected components of the induced subgraph of $T(G)$ is at most $12.$
\end{conjecture}


\subsection{Semiplanar graphs with boundary}\label{s:boundary}

In this section, we study the total curvature for semiplanar graphs with boundary. 
Let $(V,E)$ be a graph topologically embedded into a surface $S$ with boundary. We call $G=(V,E,F)$ the semiplanar graph with boundary, where $F$ is the set of faces induced by the embedding. Here we consider semiplanar graphs which satisfy tessellation properties, i.e. (i), (iii) in the introduction, and the condition (ii) replaced by the following
\begin{enumerate}[(ii')] \item The boundary of $S,$ denoted by $\partial S$, consists of edges of the graph. Each edge in the graph, by removing two end-vertices, is contained either in $\partial S$ or in the interior of $S$, denoted by $\mathrm{int}(S).$ It is incident to one face in the first case and to two different faces in the second.
\end{enumerate} Moreover, we always assume that $2\leq \deg(x)<\infty$ for any $x\in V\cap \partial S,$
$3\leq \deg(x)<\infty$ for any $x\in V\cap \mathrm{int}(S),$ and $3\leq \deg(\sigma)<\infty,$ for $\sigma\in F.$ Let $S(G)$ denote the polygonal surface for a semiplanar graph $G$ embedded into $S$ with boundary.
We define the combinatorial curvature as follows:
\begin{itemize}\item For any vertex $x\in \mathrm{int}(S),$ $\Phi(x)$ is defined as in \eqref{def:comb}.
\item  For any vertex $x$ on $\partial S$, 
\begin{equation}\label{def:curvbound}\Phi(x)=1-\frac{\deg(x)}{2}+\sum_{\sigma\in F:x\in \overline{\sigma}}\frac{1}{\mathrm{deg}(\sigma)}.\end{equation} 
\end{itemize} It is easy to see that for any $x\in V\cap \partial S,$
\begin{equation}\label{eq:gaubonboundary}\Phi(x)=\frac{\pi-\theta_x}{2\pi},\end{equation}
where $\theta_x$ is the inner angle at $x$ w.r.t. $S$. 

The patterns of a vertex on $\partial S$ with nonnegative curvature are given in the following proposition.
\begin{prop}\label{prop:curvatureboudary} Let $G$ be a semiplanar graph with boundary. If a vertex $x\in\partial S$ such that $\deg(x)\geq3$ and $\Phi(x)\geq0$, then the possible patterns of $x$ are
\begin{equation}\label{dx=2}
(3,3), (3,4), (3,5), (3,6), (4,4), (3,3,3),
\end{equation}
where the corresponding curvatures are
\begin{equation}\label{c_bd}
\frac{1}{6}, \frac{1}{12}, \frac{1}{30}, 0,0,0.
\end{equation}
\end{prop}

One can show that a semiplanar graph with boundary $G$ has nonnegative combinatorial curvature if and only if the polygonal surface $S(G)$ has nonnegative sectional curvature in the sense of Alexandrov, see \cite{BuragoGromovPerelman92,BuragoBuragoIvanov01,HJL15}. For a semiplanar graph $G$ with boundary and with nonnegative combinatorial curvature, we consider the doubling constructions of $S(G)$ and $G,$ see e.g. \cite{Perelman91}. Let $\widetilde{S(G)}$ be the
double of ${S}(G)$, that is, $\widetilde{S(G)}$ consists of two copies of ${S}(G)$ glued along the boundary $\partial S(G)$, which induces the doubling graph of $G,$ denoted by $\widetilde{G}.$ Note that $\widetilde{S(G)}$ is isometric to $S(\widetilde{G})$. One can prove that $\widetilde{G}$ has nonnegative curvature, which can be derived from either the definition of the curvature on the boundary \eqref{def:curvbound} or Perelman's doubling theorem for Alexandrov spaces \cite{Perelman91}, and $\Phi(\widetilde{G})=2\Phi(G).$
This yields that $\widetilde{G}$ is a planar graph with nonnegative combinatorial curvature if the total curvature of $G$ is positive, which yields that $S(G)$ is homeomorphic to a half-plane with boundary. For our purposes, we always assume that $S$ is homeomorphic to a half-plane with boundary. Note that the double graph $\widetilde{G}$ may have vertices of degree two if $G$ has some vertices of degree two on $\partial S.$ 
In the following lemma, we deal with the case that there is a vertex on the boundary of degree two. \begin{lemma}[Lemma~5.2 in \cite{HuaSu17a}]\label{lem:deg2}
Let $G$ be an infinite semiplanar graph with boundary and with nonnegative curvature, and $x$ be a vertex on $\partial S$ with $\deg(x)=2.$ Then the face, which $x$ is incident to, is of degree at most $6$ and $$\Phi(x)\geq\frac{1}{6}.$$ \end{lemma}

By Theorem~\ref{mainthm} and Theorem~\ref{mainthm2}, we prove that the first gap of the total curvature of semiplanar graphs with boundary and with nonnegative curvature is $\frac{1}{12}.$
\begin{corollary}\label{coro:withboundary}
Let $G$ be an infinite semiplanar graph with boundary and with nonnegative curvature. Then the total curvature of $G$ is at least $\frac{1}{12}$ if it is positive. Moreover, there is a semiplanar graph with boundary and with nonnegative curvature whose total curvature is $\frac{1}{12}.$
\end{corollary}


\begin{figure}[htbp]
\begin{center}
\includegraphics[width=0.6\linewidth]{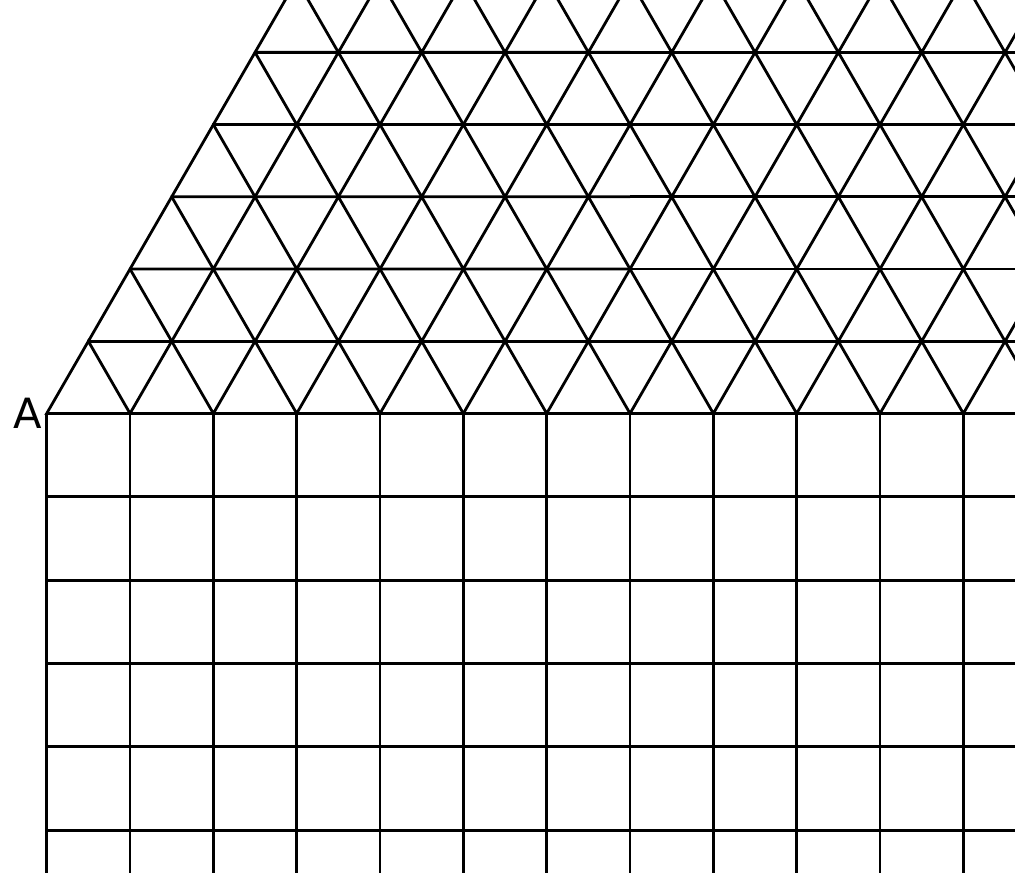}
\caption{\small A graph with boundary and with total curvature $\frac{1}{12}$.}
\label{fig-half1}
\end{center}
\end{figure}

\begin{proof} On one hand, we prove that $\Phi(G)\geq \frac{1}{12}$ if $\Phi(G)>0.$ By Lemma~\ref{lem:deg2}, it suffices to consider that $\deg(x)\geq3,$ for all $x\in V.$ Otherwise, we already have $\Phi(G)\geq \frac16.$

Consider the doubling constructions of $S(G)$ and $G,$ $\widetilde{S(G)}$ and $\widetilde{G}$. 
Note that $\Phi(G)=\frac12\Phi(\widetilde{G})$ and $\widetilde{G}$ is an infinite planar graph with nonnegative curvature since $\Phi(G)>0.$ By Theorem \ref{mainthm2}, the possible values of total curvature of $\widetilde{G}$ are $\{\frac{i}{12}|i\in\mathbb{Z}, 0\leq i\leq12\}.$ We only need to exclude the case of $\Phi(\widetilde{G})=\frac{1}{12}$. 



Suppose that $\Phi(\widetilde{G})=\frac{1}{12},$ then by the rigidity part of Theorem \ref{mainthm}, we have two cases as follows:
\begin{enumerate}
\item $S(\widetilde{G})$ is isometric to a cone with the apex angle $\theta=2\arcsin\frac{11}{12}$. In this case, there is only one vertex $x\in\widetilde{G}$ with non-vanishing curvature $\frac{1}{12}$. By the symmetry, $x\in\partial S$ and $\Phi(x)=\frac{1}{24}$ in $G$. This is impossible by \eqref{c_bd} in Proposition~\ref{prop:curvatureboudary}.
\item $S(\widetilde{G})$ is isometric to a ``frustum" with a hendecagon base. In this case, there exist eleven vertices on a hendecagon with the curvature $\frac{1}{132}$. These vertices cannot appear on $\partial S$ by \eqref{c_bd} in Proposition~\ref{prop:curvatureboudary}. Hence they are contained in $\mathrm{int}(S)$, i.e. on one side of $\widetilde{S(G)}$, which yields $\Phi(G)\geq\frac{1}{12}$ and $\Phi(\widetilde{G})\geq\frac{1}{6}$. This is a contradiction.
\end{enumerate}
So that $\Phi(\widetilde{G})>\frac{1}{12}$ which implies $\Phi(\widetilde{G})\geq\frac{2}{12}.$ Hence $\Phi(G)=\frac12\Phi(\widetilde{G})\geq \frac{1}{12}.$ 

On the other hand, we may construct a graph $G$ with boundary and total curvature $\frac{1}{12},$ see Figure \ref{fig-half1}. 
\end{proof}

The classification of the graph/tessellation structures of semiplanar graphs with boundary and with nonnegative curvature whose total curvatures attain the first gap of the total curvature was obtained in \cite{HuaSu17a}.







\bigskip
{\bf Acknowledgements.} We thank Tam\'as R\'eti for many helpful discussions on the problems.

B. H. is supported by NSFC, grant no. 11401106. Y. S. is supported by NSF of Fujian Province through Grants 2017J01556, 2016J01013, JAT160076.

\bibliography{Reti-ref-sec}
\bibliographystyle{alpha}

\end{document}